\documentclass[12pt]{amsart}
\usepackage[francais]{babel}
\usepackage{graphics}
\usepackage{graphicx}
\usepackage{pdfpages}

\title{Sur un article de 1954 sign\'e N. Cuesta\\ \tiny une traduction} 
\author{Labib Haddad}
\address{120 rue de Charonne, 75011 Paris, France}
\email{labib.haddad@wanadoo.fr}

\usepackage{amssymb}
\usepackage{amsmath}
\usepackage{endnotes}

\newcommand{\su}{\subsection*}
\newcommand{\head}{\section*}
\newcommand{\noi}{\noindent}

\newcommand{\Ž}{\'e}

\newcommand{\ˆ}{\`a}
\newcommand{\}{\`u}

\newcommand{\ž}{\^u}

\newcommand{\f}{\varphi}

\newcommand{\leqs}{\leqslant}

\newcommand{\geqs}{\geqslant}

\newcommand{\guil}{\guillemotleft}  
\newcommand {\guir}{\guillemotright}

\newcommand{\ali} {\begin{aligned}}   
\newcommand{\ala} {\end{aligned}}

\newcommand {\et}{\ \text{et}\ }

\newcommand {\ou}{\ \text{ou}\ }

\newcommand {\si}{\ \text{si}\ }

\newcommand {\pour}{\ \text{pour}\ }

\newcommand {\pourtout}{\ \text{pour tout}\ }

\newcommand{\inc}{\subset}

\newcommand{\inq}{\subseteq}

\newcommand{\bc}{\begin{cases}}
\newcommand{\ec}{\end{cases}}

\newcommand{\ba}{\begin{array}}
\newcommand{\ea}{\end{array}}

\newcommand{\uds}{\underset{*}}
\newcommand{\ovs}{\overset{*} }
\newcommand{\udl}{\underline}
\newcommand{\ovl}{\overline}

\begin{document}
\maketitle
\thispagestyle{empty}
\markboth{}{}

\

\

\head{Preamble}

\

\

Here is a translation from Spanish to French of a paper dating back to 1954, published by N. Cuesta.

\

The paper deals mainly with partially, and totally, ordered sets. Two subjects are specially dealt with: Construction of new ordered sets starting from a family of those. Completion of ordered sets by tools akin to Dedekind cuts. Curiously enough, the so-called surreal numbers (later defined by Conway, in 1974) are already there, thirty years before.

\

{\bf The translation is done for research purposes. It is not intended to obtain financial gains}. I tried, with the precious help of my grand niece, Aya Nay Haddad, to find the \lq\lq copyright owner", without success. The librarians at Columbia University gave her some advice about the matter.

\

So, to my knowledge there is no copyright owner. The author of the paper, Norberto Cuesta, died in 1989 (last century). 

\

Just in case there still are copyright owners, I beg them to contact me on my email address.

\

Tous mes remerciements \ˆ Aya Nay, pour son aide, sans qui je n'aurais jamais r\Žussi \ˆ d\Žm\ler les probl\mes de droit d'auteur ni p\Žn\Žtrer les arcanes de ce domaine complexe !

\newpage

\head{Introduction}

\

\

Il s'agit d'une traduction en fran\c cais de l'une des premi\res publications de Norberto Cuesta Dutari. C'est un article, \Žcrit en espagnol, paru dans la Revista de la Real Academia de Ciencias Exactas, Fisicas et Naturales  de Madrid, intitul\Ž {\it  Algebra Ordinal}, publi\Ž en 1954 et sign\Ž, simplement, {\sc N. Cuesta} [ 3 ]. Il est tir\Ž de la th\se qu'avait soutenue le jeune math\Žmaticien en 1943. 

\

J'en ai eu connaissance en parcourant le livre tr\s document\Ž et pr\Žcieux de {\sc Alling} [ 1 ] : comme son titre  le pr\Žcise,  ce livre est un trait\Ž qui porte sur les nombres surr\Žels. Mais le titre ne le dit pas, c'est ainsi que l'on appelait couramment les nombres introduits par {\sc Conway} dans les ann\Žs 1970. 

\

{\sc Alling} donne la r\Žf\Žrence \ˆ  l'article de {\sc Cuesta} mais ne signale point o\ on peut le trouver. 

\

Aucune des biblioth\ques que j'ai contact\Žes ne poss\dait ce num\Žro particulier de la Revista. Il faut rappeler qu'en ce temps-l\ˆ l'Espagne vivait toujours quelque peu isol\Že du reste du monde. Trouver une copie  de l'article s'est r\Žv\Žl\Ž \tre une t\‰che tr\s difficile. Cette traduction est quasiment une affaire de famille. 

\

En effet, n'arrivant pas \ˆ obtenir une copie de l'article de Cuesta, je me suis adress\Ž \ˆ mon ami et excellent coll\gue, Charles Helou, Professeur de math\Žmatiques \ˆ l'Universit\Ž d'\'Etat de Pensyvalnie. Il  est  arriv\Ž, en sollicitant les biblioth\caires de son institution  \ˆ obtenir enfin une copie du texte, si difficile \ˆ d\Žnicher. Comme on peut le voir sur la couverture de la revue, reproduite ci-dessous, l'exemplaire  provenait de la biblioth\que du D\Žpartement de l'Agriculture des \' Etats-Unis, curieusement. Comme on peut le remarquer \Žgalement, en dessous, la copie,  avec les courbures de ses pages, n'est pas facile \ˆ exploiter, sur un ordinateur.

\includepdf{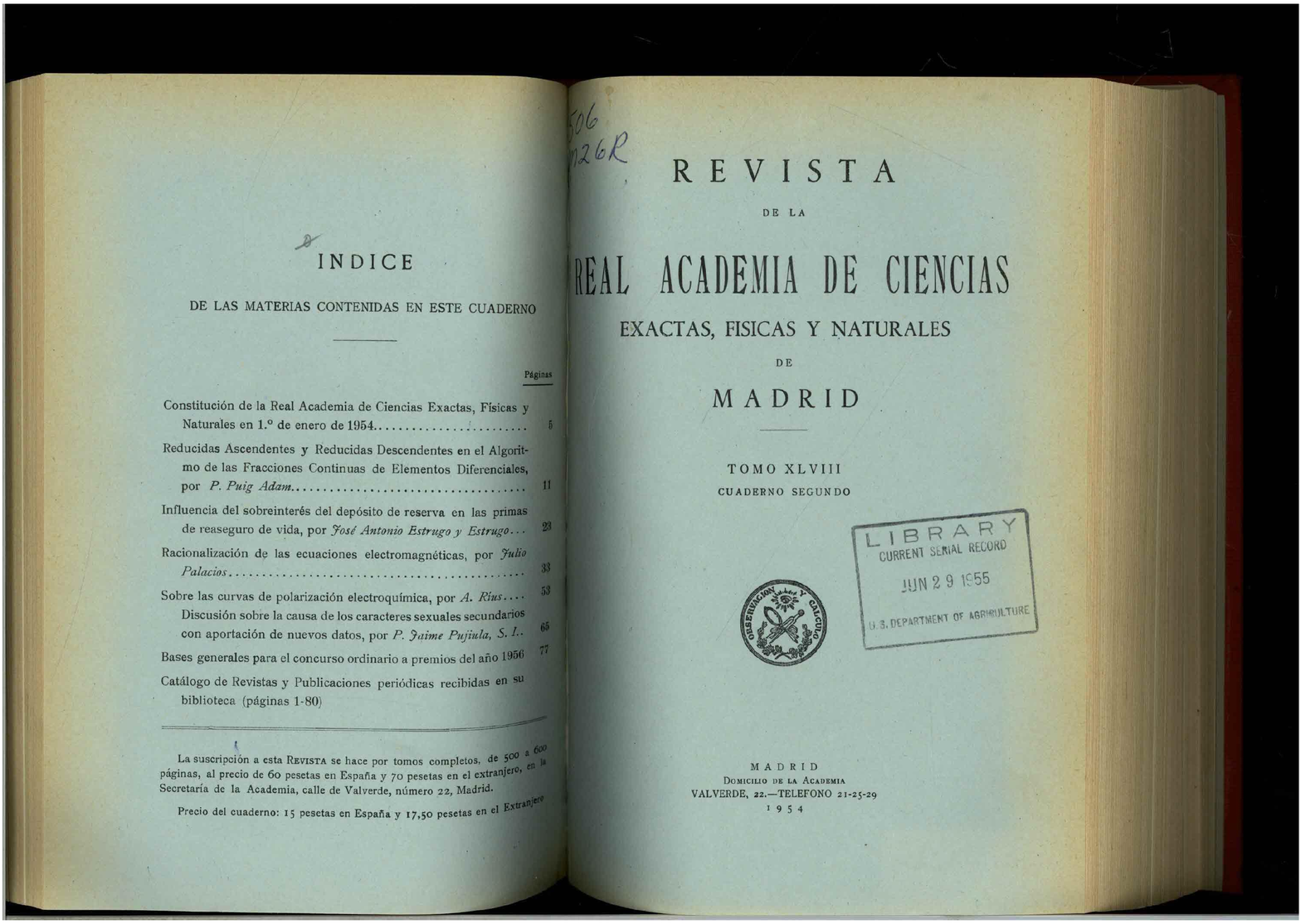}

\includepdf{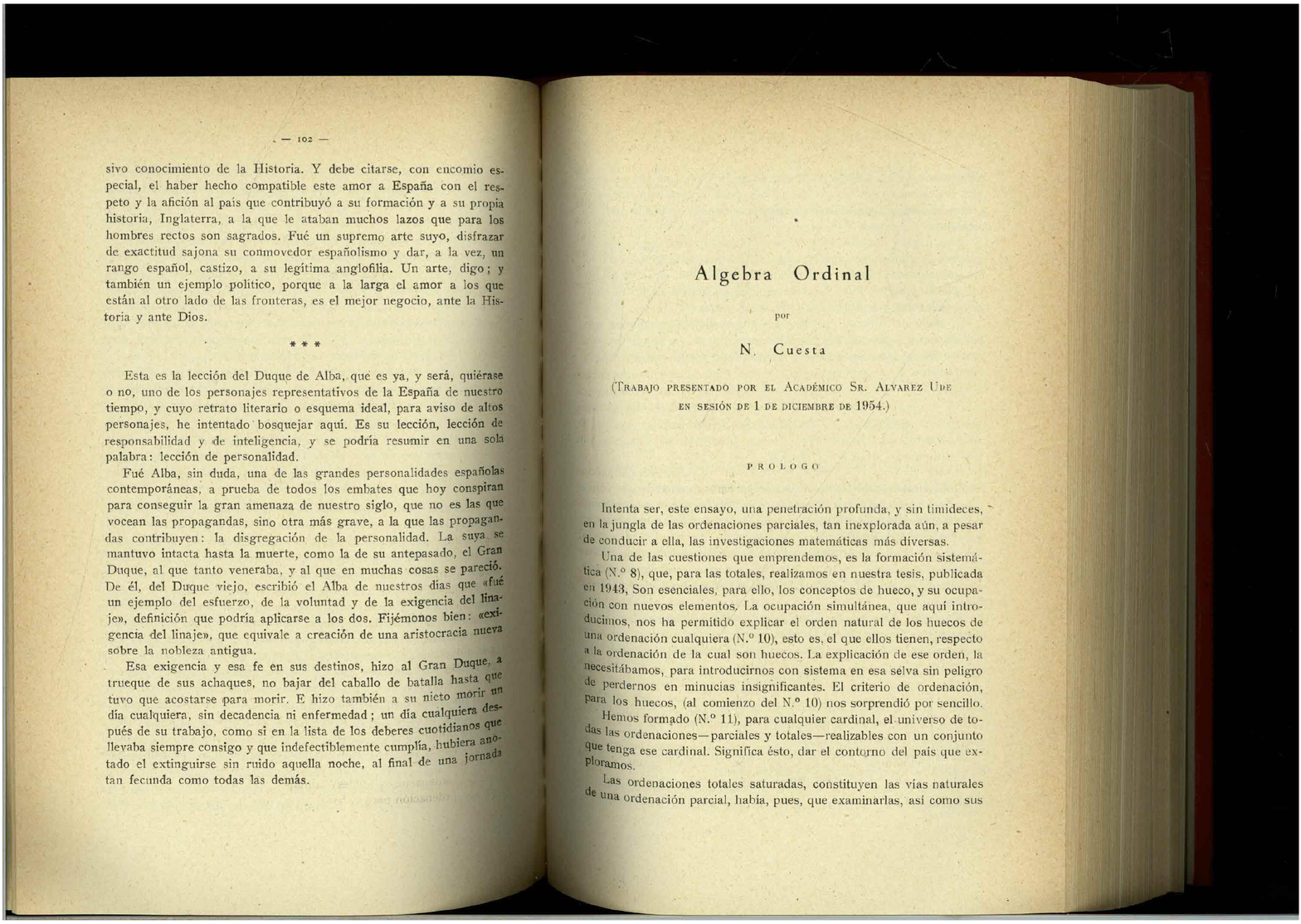}

\

Aussi, est-ce ma ni\ce, Samar Haddad, qui a r\Žussi \ˆ l'aide d'outils informatiques sp\Žciaux \ˆ faire appara\"tre le texte aplani et bien plus maniable, comme on peut le constater sur la reproduction ci-dessous.

\

Cela ne r\Žsolvait pas, pour autant, tous les probl\mes soulev\Žs par la traduction. Il fallait encore retrouver un peu de la {\bf saveur} particuli\re de la langue espagnole de Cuesta. C'est ma ch\re \Žpouse, Claude Boisnard Haddad, qui m'a amplement aid\Ž \ˆ tenter de la reproduire.

\

On n'a  pas corrig\Ž, dans les formules, les quelques coquilles et les rares omissions. On a essay\Ž de les reproduire telles quelles, autant que faire se peut,  afin de garder son authenticit\Ž au texte. Elles seront ais\Žment corrig\Žes par le lecteur.

\

J'ai essay\Ž de respecter le plus scrupuleusement possible les notations de l'original et sa mise en page. Tous les d\Žfauts et les erreurs qui y persistent sont de mon fait. Je ne m'en vante pas mais les revendique !

\

\su{Remerciements} Je remercie chaleureusement Charles Helou, pour sa pers\Žv\Žrance et sa perspicacit\Ž. Toute ma reconnaissance \ˆ ma ni\ce, Samar Haddad, \ˆ qui ma gratitude est enti\rement acquise, qui a consacr\Ž beaucoup de son temps pr\Žcieux et son savoir-faire pour arriver au r\Žsultat voulu. Elle m'a \Žgalement donn\Ž de tr\s bons conseils au sujet de l'usage du format pdf pour la reproduction des figures. 

Enfin, et c'est essentiel, \ˆ mon \Žpouse, si d\Žvou\Že, comp\Žtente et patiente, j'adresse un immense merci !

\

\
 
 \head{Bibliographie}
 
\

\noi [1] {\sc Norman L. ALLING}, {Foundations of analysis over surreal number fields}, Mathematics Studies 141, xvi + 373 pp., North Holland, 1987.

\

\noi [2] {\sc J. H. CONWAY}, {On numbers and games}, ix + 230 pp.,
Academic Press Inc. 1976, reprinted 1979. 

\

\noi [3] {\sc N. CUESTA}, {\sl Algebra ordinal}, Revista de la Real Academia De Ciencias Exactas, Fisicas y Naturales, {\bf 58} \no 2 (1954) 103-145.

\includepdf{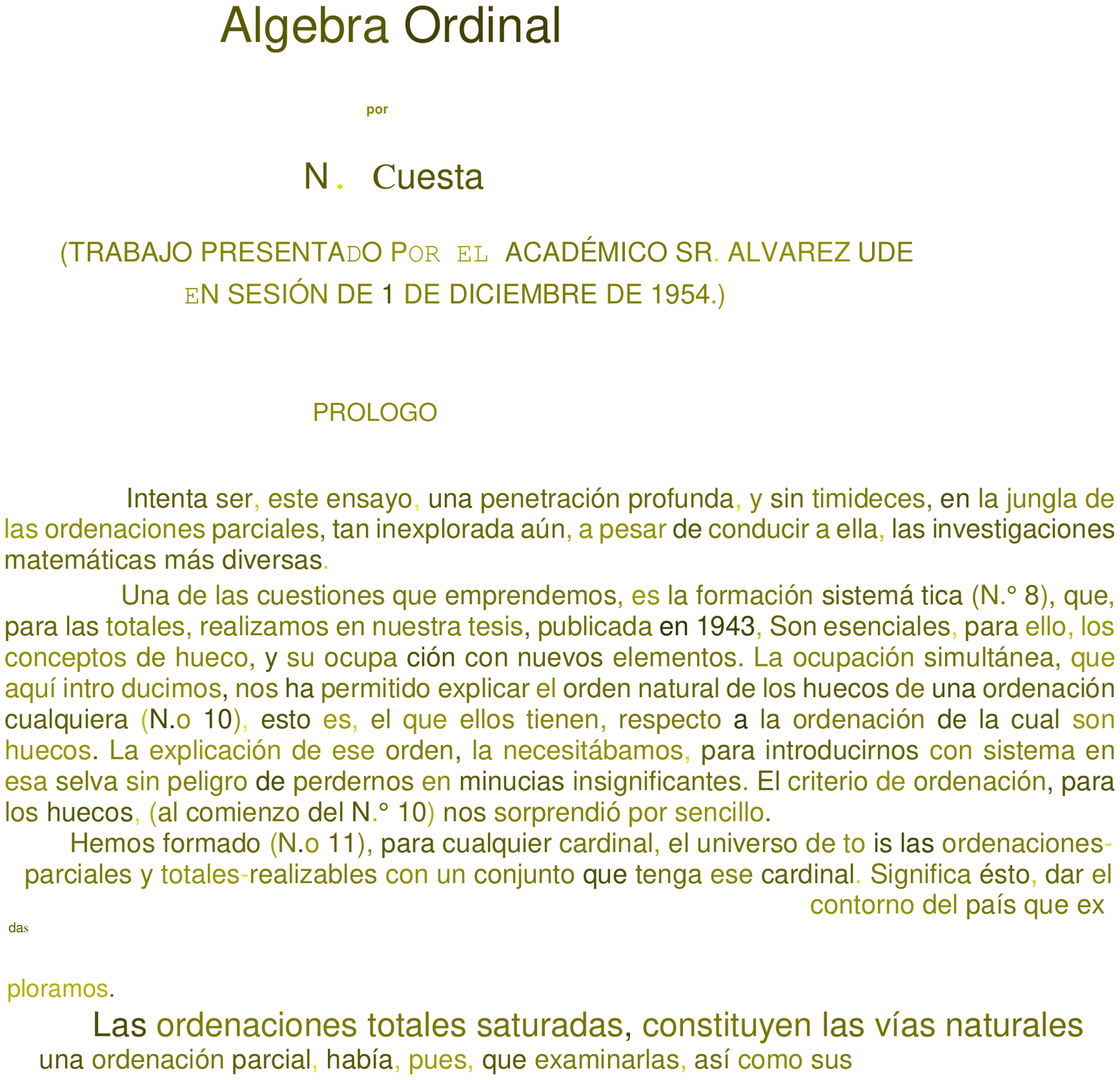}

\

\

\centerline{\bf  UNE TRADUCTION}

\

\head{Alg\bre ordinale}

\head{par}

\head{N. Cuesta}

\head {(Texte pr\Žsent\Ž par l'Acad\Žmicien M. \'ALVAREZ UDE
\ˆ la session du 1er d\Žcembre 1954)}

\

\head{Prologue}

\

\

Cet essai se pr\Žsente comme une incursion profonde, et audacieuse, dans la jungle des ordres partiels,  jusqu'\ˆ pr\Žsent tellement inexplor\Že, bien que les recherches math\Žmatiques les plus diverses y conduisent.

L'une des questions que nous abordons est la formation syst\Žmatique (\no 8) que, pour l'ensemble, nous avons r\Žalis\Že dans notre th\se, publi\Že en 1943. Sont essentielles pour cela les notions de trou et de remplissage d'un trou par de nouveaux \Žl\Žments. Le remplissage simultan\Ž que nous introduisons ici, nous a permis d'expliquer l'ordre naturel des trous d'un ordre quelconque (\no  10), c'est-\ˆ-dire celui qu'ils ont par rapport \ˆ l'ordre dont ils sont les trous. 
Pour trouver l'explication de cet ordre, il nous fallait entrer syst\Žmatiquement dans cette jungle sans risquer de nous perdre dans des d\Žtails insignifiants. Le crit\re de l'ordre, pour les trous, (voir le d\Žbut du \no 10) nous a surpris par sa simplicit\Ž.

Nous avons construit (\no 11), pour tout cardinal, l'univers de tous les ordres -- partiels et totaux -- r\Žalisables sur un ensemble  ayant ce cardinal. Cela signifie, pour nous, avoir cern\Ž la contr\Že que nous avons explor\Že.

Les ordres totaux satur\Žs \Žtant les voies naturelles vers un ordre partiel, il fallait les examiner, ainsi que leur
croisement et les trous par lesquels ils passent (\no 12). Nous \Žtudions \Žgalement les transversales compl\tes (\no 13) car, si ce sont des ensembles amorphes, au fond, ils entretiennent des relations ordinales int\Žressantes avec leur alentour.

Pour notre \Žtude, nous avons eu recours \ˆ diverses figures que nous retenons pour leur grand pouvoir suggestif que nos lecteurs appr\Žcieront bien que celles que l'on peut dessiner ne couvrent qu'une infime partie des ordres possibles. Ces figures incitent \ˆ  \Žtudier la gen\se des ordres partiels par fusion conjointe d'ordres totaux (\no 14).

Le titre de l'article est d\ž aux op\Žrations g\Žn\Žratrices que nous introduisons, et qui constituent le langage avec lequel l'appr\Žhension mentale de ces objets subtils et insaisissables peut \tre r\Žalis\Že rapidement et en toute s\Žcurit\Ž.

\

\head{CHAPITRE 1 LES STRUCTURES BINAIRES}
				
\

1.--{\sc Les structures binaires en g\Žn\Žral}

\

D\Žsignant par $M \times P$ l'ensemble des paires ordonn\Žes $mp$ form\Žes par les \Žl\Žments g\Žn\Žriques respectifs $m$ et $p$ de ces ensembles, $M \times M$ repr\Žsentera les paires ordonn\Žes form\Žes \ˆ l'aide du seul ensemble $M$.

Sur un ensemble $M$, la structure $(M\&)$ consistera en  la donn\Že, \ˆ l'aide du signe de relation  \&, lorsque, pour certaines paires ordonn\Žes appartenant \ˆ $M \times M$, le signe de relation peut s'intercaler,  en \Žcrivant $a \ \&\ b$.

Les \Žl\Žments {\it r\Žflexifs} de $(M\&)$ sont ceux pour lesquels on a $r \ \&\ r$. Nous dirons que les autres sont {\it irr\Žflexifs}.

Lorsque deux \Žl\Žments $a$ et $b$ v\Žrifient simultan\Žment
$$a \ \&\ b \et b \ \&\ a$$
nous dirons que ces deux  \Žl\Žments forment une paire {\it sym\Žtrique} de la structure binaire $(M \&)$. En particulier, on compte parmi eux les paires form\Žes d'\Žl\Žments r\Žflexifs.

Si de ces deux relations une seule est satisfaite, on dira que la paire form\Že des deux \Žl\Žments $a$ et $b$ est {\it asym\Žtrique}.

Si on ne pouvait \Žcrire aucune, nous dirions que les deux \Žl\Žments $a$ et $b$ sont {\it incomparables}. En particulier, un \Žl\Žment irr\Žflexif est incomparable \ˆ lui-m\me. Dans les autres cas, nous dirons que les \Žl\Žments $a$ et $b$ sont {\it comparables}.

Pour d\Žsigner diverses structures binaires, surtout lorsque l'ensemble structur\Ž est le m\me, nous mettrons des accents, ou des indices, au signe de la relation.

Nous dirons que la structure binaire $(M \&_1)$ est {\it plus faible} que  la structure $(P \&_2)$ lorsque :

1.$^\circ$) $M$ est un sous-ensemble de $P$.

2.$^\circ$) deux \Žl\Žments $a$ et $b$ de $M$ qui v\Žrifient $a \ \&_1\ b$ v\Žrifient \Žgalement $a \ \&_2 \ b$. Nous \Žcrirons  $(M \&_1) \leqs (P \&_2)$.

Une paire de sous-ensembles $A$ et $B$ non vides, sans \Žl\Žments communs, sont dits {\it non connect\Žs} dans $(M \&)$ lorsque, quels que soient $a$ et $b$, respectivement de $A$ et $B$, sont incomparables.

On dira que la structure $(M \&)$ est {\it connexe} lorsqu'il n'existe aucun sous-ensemble $X$ de $M$  non connect\Ž \ˆ son compl\Žmentaire $M - X$.

{\it R\Žunir} les structure binaires $(M^j \&)_{j\in J}$ c'est former la nouvelle structure $(\sum_{j\in J} M^j \&)$ sur l'ensemble r\Žunion o\ l'on aura $a \ \& \ b$ lorsque l'on avait $a \ \& \ b$ dans l'une au moins des structures $(M^j \&)$.

Toute structure binaire $(M \&)$ est r\Žunion d'autres structures $(M^j \&)$ dont chacune est connexe et qui sont non connect\Žes deux \ˆ deux.

En effet : d\Žsignons par $a\f$ l'ensemble obtenu en adjoignant \ˆ l'\Žl\Ž\-ment $a$ tous ceux qui lui sont comparables. En r\Žit\Žrant l'op\Žration $\f$, on obtient la suite
$$a  \inq a \f_1 \inq a \f_2 \inq a \f_3 \inq  \dots \inq a \f_n \inq \cdots$$
o\ 
$$a\f_n \equiv (a\f_{n-1})\f$$
\Žtant sous-entendu que
$$a\f_1 \equiv a\f$$
On \Žcrit 
$$M^1 \equiv \lim_{n\to \omega_0} (a\f_n)$$

$M^1$ est connexe : en effet; d\Žsignons, de deux sous-ensembles compl\Ž\-mentaires, par $X$ celui qui contient $a$; par $Y$ son compl\Ž\-mentaire. Si ces deux sous-ensembles \Žtaient non connect\Žs, on aurait, \Žtant entendu que
$$a\f_0 \equiv a$$
pour tout entier $i\geqs 0$ :
\noi si
$$a\f_i \inq X \ , \ a\f_{i+1}. Y \equiv {\rm O}  \  \therefore \ a\f_{i+1} \inq X$$
donc
$$M^1 \inq X$$
et $Y$ serait vide.

$M^1$ est relativement ferm\Že; cela veut dire
$$M^1\f  \equiv M^1$$

 En effet : quel que soit l'\Žl\Žment  $m$ de $M^1$, il appara\"tra \ˆ une \Žtape $a\f_n$; de sorte que tous les $m\f$ figurent dans $a\f_{n+1}$; et puisque
$$M^1\f \equiv \sum_{m\in M^1} m\f \inq M^1$$
\Žtant \Žvident que
$$M^1 \inq M^1\f$$
des deux r\Žsulte l'identit\Ž annonc\Že.

De cette identit\Ž s'ensuit clairement  que si $M - M^1 \not\equiv 0$, $M^1$ et son compl\Žmentaire seraient non connect\Žs.

Prenant $b$ dans $M - M^1$ s'il n'est pas vide, on obtient de mani\re analogue $(M^2 \&)$ \Žgalement connexe, et ainsi de suite.

On voit ais\Žment que cette d\Žcomposition de $(M \&)$, en structures binaires connexes, deux \ˆ deux non connect\Žes, est univoque.

Nous \Žcrirons $A  \ \& \ B$ \Žtant donn\Žs $A$ et $B$  sous-ensembles  de $M$, lorsque, quels que soient $a$ et $b$, \Žl\Žments respectivement de $A$ et $B$, $ a \ \&\ b$ est vrai. Soit $M_1$ le syst\me des sous-ensembles de $M$, la structure $(M \&)$ nous permet de d\Žfinir $(M_1 \&)$. \`A cette op\Žration g\Žn\Žratrice nous donnerons le nom de {\it sousconjoindre} dans $(M <)$. \'Evidemment, on peut r\Žit\Žrer et obtenir $(M_2 \&)$. Si $M_{\omega_0}$  d\Žsignait l'ensemble qui comprend tous les pr\Žc\Ždents, $(M_{\omega_0} \&)$ d\Žsignerait  n\Žcessairement d\Žfinie sur $M_{\omega_0}$,  la relation $\&$ suivant les \Žtapes pr\Žc\Ždentes; donc
$$(M_{\omega_0} \&) \equiv\ \left(\sum_{O\leqs n < \omega_0} M_n \&\right)$$

Cela pourrait \tre poursuivi encore, donnant un sens au symbole $(M_\alpha \&)$ pour n'importe quel ordinal $\alpha$ donn\Ž.

Nous appellerons {\it ligne compl\te} de $(M \&)$ la sous-structure 
$(L \&)$ induite sur un sous-ensemble $L$ de $M$ lorsqu'il n'existe pas dans $L$ de  paires d'\Žl\Žments incomparables et que chaque \Žl\Žment de $M-L$ est incomparable \ˆ chacun des \Žl\Žments de $L$.

Nous appellerons {\it transversale compl\te} de $(M \&)$ tout sous-ensemble $T$ de $M$ dont les paires d'\Žl\Žments diff\Žrents sont incomparables, chaque \Žl\Žment de $M-T$ \Žtant comparable \ˆ certains \Žl\Žments de $T$.

Les \Žl\Žments de $M$ solutions de la relation
$$x \ \& \ a$$
pour un $a$ donn\Ž de $M$, forment un ensemble que nous d\Žsignerons $\uds a$. De mani\re analogue, \ˆ droite, on aura le symbole $\ovs a$.

{\it Fermer \ˆ gauche} un ensemble $A$ voudra dire lui adjoindre les \Žl\Žments de $\uds a$ pour chacun de ses \Žl\Žments $a$. En le d\Žsignant par $\udl A$, on aura
$$\udl A \equiv A + \sum_{a\in A} \uds a$$

De m\me, on d\Žfinira $\ovl A$ et on aura
$$\ovl A \equiv A + \sum_{a\in A} \ovs a$$

{\it Inverser} la structure $(M \&)$ c'est construire la nouvelle structure $(M \&')$ o\ $a \ \&' \ b$ \Žquivaut \ˆ $b \ \& \ a$.

Il convient d'observer que, en donnant le nom $P$ au sous-ensemble de $M\times M$ form\Ž des paires li\Žes par le symbole de relation $\&$, et $N$ \Žtant l'ensemble des \Žl\Žments qui appara\"ssent dans les paires de $P$, le seul r\Žellement structur\Ž est l'ensemble $N$. En passant de $(N \&)$ \ˆ $(M \&)$, les \Žl\Žments de $M-N$ n'ont aucune relation entre eux, ni avec les \Žl\Žments de $N$. On peut dire que $(N \&)$ \Žtait le {\it noyau structur\Ž} tandis que $M-N$ \Žtait un {\it r\Žsidu amorphe}.

Par cons\Žquent, il sera plausible de consid\Žrer comme {\it isomorphes au sens large }deux structures $(M \&_1)$ et $(Q \&_2)$ lorsque les structures de leurs noyaux le sont au sens strict.

\

2.--{\sc Les structures binaires transitives}

\

\`A l'aide du symbole $(M \to)$ nous repr\Žsentons les structures binaires transitives, autrement dit, celles pour lesquelles
$$a \to b \et b\to c \ \text{impliquent} \ a \to c $$

Un exemple suggestif de ces structures est fourni par la figure jointe

\

\includegraphics[width=4.8in]{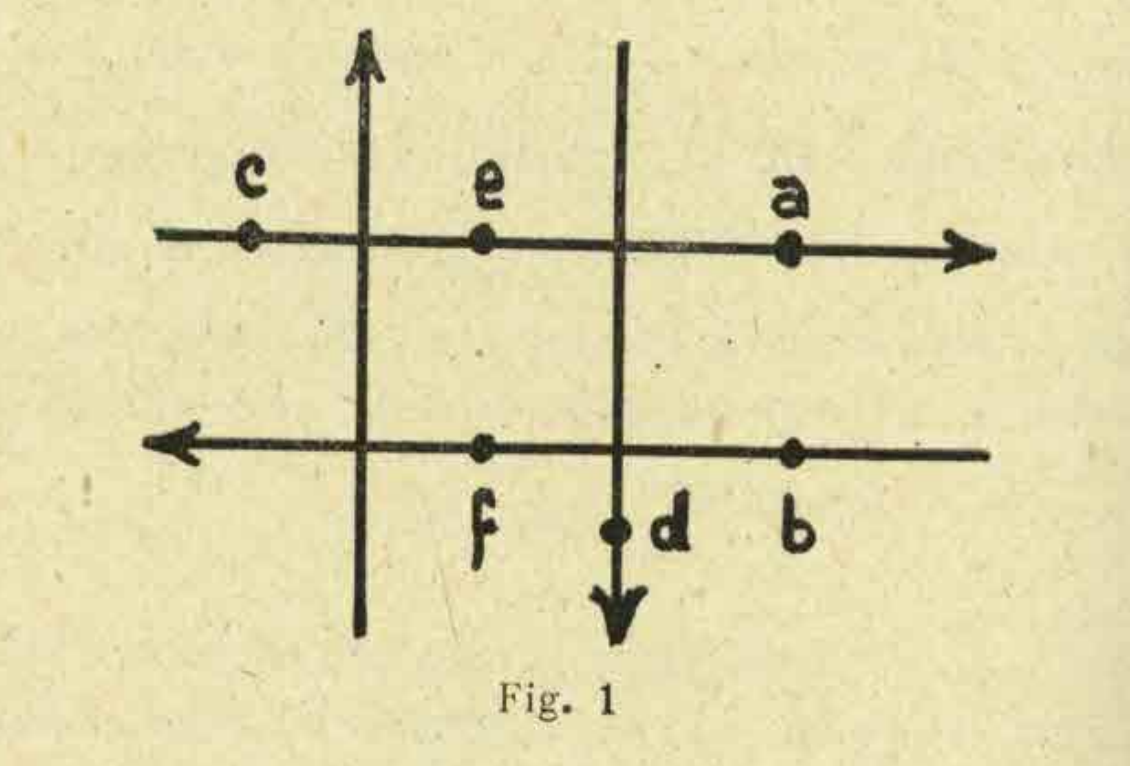}

\noi $M$ est constitu\Ž des points qui appartiennent aux droites de la figure. Pour deux d'entre eux, nous poserons $x\to z$, lorsque l'on peut aller du premier au second, en faisant un parcours effectif, jamais \ˆ contre-courant des fl\ches des
lignes droites. Sur cette figure, $a b$ constitue une paire incomparable, $c d$ une paire comparable, donc
$c\to d $, mais asym\Žtrique, tandis que $e f$ constitue une paire comparable sym\Žtrique  donc ils v\Žrifient
$$e \to f \ \et \ f \to e$$
Les \Žl\Žments $a, b, c, d,$ sont eux tous irr\Žflexifs. Les \Žl\Žments $e, f,$ sont tous deux r\Žflexifs.

En intercalant le symbole de la relation  $\to$ dans les diverses paires de $M\times M$, on obtient des relations ou {\it propositions}.

Nous dirons que l'une de ces propositions est {\it vraie} lorsque la relation correspondante est v\Žrifi\Že dans $(M \to)$; sinon, nous dirons qu'elle est {\it fausse}. Ainsi, la structure $(M \to)$ constitue le {\it crit\re de v\Žrit\Ž}.

La loi transitive fonctionne comme un {\it m\Žcanisme d\Žductif} pour le syst\me de toutes les propositions. Lorsqu'on l'applique \ˆ des propositions vraies, il nous en donne d'autres vraies.

Quel que soit le sous-ensemble de $M \times M$, en intercalant le symbole de la relation $\to$ dans ses paires, nous obtenons un sous-ensemble $P$ de l'ensemble des propositions. En appliquant directement la loi de transitivit\Ž \ˆ ses paires de propositions, \ˆ l'aide d'un \Žl\Žment moyen commun, on obtient un nouvel ensemble $P \f$. En r\Žit\Žrant le proc\Žd\Ž sur cet ensemble de propositions et ainsi, successivement, on aura
$$P \inq P\f _1\inq P \f_2 \inq P \f_3 \inq \dots \inq P \f_n \inq \cdots$$
o\, pour pour chaque ordinal fini $n$, on interpr\te
$$P\f_n \equiv (P\f_{n-1})\f \ \text{sous-entendu} \ P\f_1 \equiv P\f$$

S'il se v\Žrifie que
$$P\f_{n-1} \equiv P\f_n$$
dans toutes les \Žtapes suivantes, on n'en obtiendra rien de nouveau; mais, si cela ne se produit \ˆ aucune des \Žtapes finies, l'ensemble limite, que nous d\Žsignerons par $P \f_{\omega_0}$, est transitivement ferm\Ž. En effet : si lui appartiennent
$$a \to b \ \et \ b \to c$$ 
toutes deux apparaissant ensemble, pour la premi\re fois, en $P \f_n$; alors
$$a \to c$$
appara\"t, si ce n'est fait avant, assur\Žment en $P\f_{n+1}$; par cons\Žquent, $a\to c$ appartient \ˆ $P_{\omega_0}$.

Quelle que soit l'\Žtape de la fermeture, comme il a \Žt\Ž dit $\leqs \omega_0$, nous d\Žsignerons par $P\ovl\f$ l'ensemble, transitivement ferm\Ž, o\ s'arr\te l'efficacit\Ž du m\Žcanisme d\Žductif. Nous appellerons cet ensemble $P\ovl\f$ la {\it fermeture d\Žductive} de $P$. \`A l'aide de cet ensemble $P\ovl\f$ on peut d\Žfinir, sur $M$, une structure $(M \to)$ dont les paires li\Žes sont pr\Žcis\Žment celles donn\Žes par $P\ovl\f$.

La subordination entre $P$ et $P\ovl\f$ d\Žtermine la structure topologique ferm\Že $(M
\times M\ovl\f)$, int\Žressante car ses sous-ensembles ferm\Žs,
d\Žtermi\-nent de mani\re biunivoque les structures binaires transitives possibles sur l'ensemble $M$.

Si la structure $(M \to)$ est donn\Že a priori, l'ensemble $P$ \Žtant inclus dans l'ensemble $W$  des propositions vraies par rapport \ˆ $(M\to)$, et si $(M \to_1)$ est la structure d\Žtermin\Že par $P\ovl\f$, la structure $(M\to_1)$ serait plus faible que $(M \to)$; en symboles
$$(M\to_1) \leqs (M \to)$$

Si elle \Žtait strictement plus faible, ce serait parce que $(W-P\ovl\f)$ ne serait pas vide. Ses  propositions vraies ne viendraient jamais par d\Žduction en partant de celles de $P$. Nous dirions que les vraies, contenues dans $(W-P\ovl\f)$  \Žtaient transcendantes par rapport \ˆ $P$.

Lorsque $P\ovl\f$ co\•ncide avec $W$, la structure engendr\Že serait celle initialement donn\Že $(M \to)$. Ayant \Žt\Ž suffisant de donner la relation $\to$ uniquement entre les paires comprises dans $P$, nous dirons que c'est une {\it base} de $(M\to)$.

En g\Žn\Žral, les sous-ensembles $X$ de $M \times M$ qui v\Žrifient 
$$X\ovl \f \equiv  W$$
sont divers. Le syst\me de leur ensemble, structur\Ž par la relation binaire d'inclusion stricte $\inc$, \Žgalement transitive, est d'un grand int\Žr\t. On dira que la base $P$ est {\it r\Žductible} lorsque l'une des ses parties strictes est aussi une base : sinon, elle sera dite {\it irr\Žductible}. Il n'est pas exclu qu'il existe plusieurs bases irr\Žductibles pour une structure binaire transitive $(M \to)$, donn\Že a priori. Lorsque, parmi l'ensemble des bases, il y en a une contenue dans chacune des autres, on dira que c'est une {\it base absolue} de $(M\to)$. A priori, il n'est pas exclu que des structures existent non seulement sans base absolue, mais m\me sans aucune base irr\Žductible. A posteriori, on sait qu'il y en a : par exemple, il est facilement prouv\Ž que l'ordre arithm\Žtique des nombres rationnels manque d'une base irr\Žductible. (*)

\noi -----------------

(*) Voir notre article  \guil Mod\les d\Žductifs ordinaux\guir \ (Rev. Mat. Hisp. Am. (IV) {\it 13} (1953).

\

\' Etant donn\Ž le syst\me de structures $(M^j \to)$, d\Žfini sur les ensembles $M^j$, d\Žtermin\Žs pour chacun des \Žl\Žments $j$ d'un certain ensemble $J$, pouvant avoir -- et c'est le cas le plus int\Žressant -- des ensembles $M^j$ avec \Žl\Žments communs, la structure
$$\left(\sum_{j\in J} M^j \to\right)$$
qui, gr\‰ce \ˆ la loi transitive, n'est pas la simple r\Žunion des $(M^j \to)$, nous l'appellerons {\it con-fusion} -- donnant \ˆ ce mot sa valeur g\Žn\Žtique -- des structures du syst\me $(M^j \to)$. Si $M'$ et $M"$ avaient des \Žl\Žments communs, au moyen de ceux-ci, et de la loi transitive, deux \Žl\Žments non communs de $M'$ et $M"$ pourraient \tre li\Žs.

\

\includegraphics[width=4.8in]{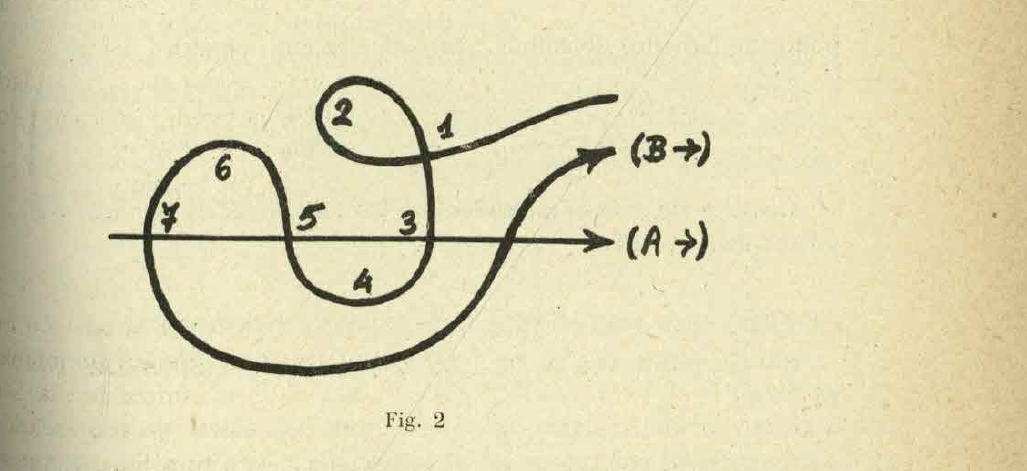}

\

Dans la figure 1, on peut consid\Žrer la structure $(M\to)$, issue de la fusion conjointe de structures semblables, situ\Žes sur chacune des droites qui la composent. Dans la figure 2 jointe, la structure $(M \to)$ est obtenue par  con-fusion des structures $(A \to)$  et $(B\to)$, situ\Žes sur chacun des traits qui la  compose. Les morceaux $(1, 2, 1)$ et $(3, 4, 5, 6, 7, 5, 3)$ contiennent les \Žl\Žments r\Žflexifs.

Les propri\Žt\Žs suivantes des structures binaires transitives m\Žritent d'\tre mentionn\Žes :

 a) Les \Žl\Žments d'une paire sym\Žtrique sont r\Žflexifs. 

En effet : $a \ \to \ b$ et  $b \ \to \ a$  impliquent aussi bien $a \ \to \ a$  que  $b\ \to \ b$. 

Par cons\Žquent, les structures binaires transitives, sans \Žl\Žments r\Žflexifs, n'ont pas de paires sym\Žtriques. Bien s\žr une structure transitive peut avoir des \Žl\Žments r\Žflexifs, sans qu'il soit n\Žcessaire qu'elle contienne des paires sym\Žtriques d'\Žl\Žments diff\Žrents.

b) Si une structure binaire est transitive, des \Žl\Žments sym\Žtriques \ˆ un troisi\me sont sym\Žtriques entre eux. La preuve est triviale.

c) Si $a \ \to \ c \ \to b$ devait \tre vraie et que les \Žl\Žments $a$ et $b$ \Žtaient sym\Žtriques, on aurait aussi $b \ \to  \ c \ \to \ a$

En effet : de $b \ \to \ a$ et $a \ \to \  c$ r\Žsulte $b \ \to \ c$
 
De m\me de $c \ \to \ b$ et $b \ \to \ a$ r\Žsulte $c\ \to a$

d) Toute ligne compl\te de la structure $(M \to)$ contient, avec l'\Žl\Ž\-ment $m$, tous ses sym\Žtriques.

En effet : soit $(L \to)$  la ligne compl\te. Les \Žl\Žments de $L$ v\Žrifient ou bien
$$m \ \to \ l'' \ou l' \ \to m$$
ou encore les deux. Soit $p$ un \Žl\Žment sym\Žtrique de $m$
$$\text{de} \ l' \ \to m \ \et m \ \to \ p \ \ \text{s'ensuit}\ l' \ \to \ p$$
$$\text{de} \ p \ \to \ m \et m \ \to \ l''\  \ \text{s'ensuit}\ p \ \to \ l''$$

Comme le voit, $p$ est comparable aux \Žl\Žments de la ligne et dans le m\me sens que $m$.

\

3.--{\sc La classification des structures binaires transitives}

\

Dans le tableau suivant, par caract\Žristiques tr\s visibles, nous classons les structures binaires transitives.

\

\noi Structures binaires transitives Symbole $(M\to)$ 

\

$\bc

1\ \text{Irr\Žflexifs} \\ 
              
\text{(pas d'\Žl\Žments  r\Žflexifs)}.\\

\text {Symbole} \ (M <).\\

\\

\text{ORDRES} \\

\\

2 \ \text{R\Žflexifs}\\

\text{(tous les \Žl\Žments sont r\Žflexifs)}\\

\text {Symbole}\ (M \leqs). \\

\\

3 \ \text{Avec \Žl\Žments r\Žflexifs et irr\Žflexifs}

\ec$

\

\

\

\noi 1\ Irr\Žflexifs $\bc  

11\ \text{Sans paires incomparables.  ({\sc ordres totaux}) }\\

\\

12 \ \text{Avec paires incomparables. ({\sc ordres partiels})} \\

\ec$

\

\

\

\noi 2 \ R\Žflexifs $\bc

21 \ \text{Sans paires asym\Žtriques. ({\sc classifications})}\\

 \text{Symbole} \ (M =)\\
 
 \\

22 \ \text{Avec paires asym\Žtriques}\bc

221 \ \text{Sans paires}\\

\text{sym\Žtriques.}\\

\\

222 \ \text{Avec paires}\\

\text{sym\Žtriques.}\\

\ec

\ec$

\

\

\

Des exemples clairs [{\sl perspicuos}] de la troisi\me classe, les figures faites en donnent, avec des \Žl\Žments r\Žflexifs et irr\Žflexifs. Cependant, cette classe a beaucoup moins d'importance que les pr\Žc\Ždentes car toute structure de cette classe devient une autre de la seconde, rien qu'en \Žcrivant
$$a\leqs b  \ \text{lorsque}\ a \ \to \ b$$
ajoutant aussi
$$a \leqs a$$
pour chaque \Žl\Žment de l'ensemble de base.

Concernant les structures de classification, nous devons en profiter pour
rectifier notre affirmation de l'article \guil Estructuras y programa de Erlangen\guir \ [{\sl Structures et programme d'Erlangen}] que, de la loi sym\Žtrique et transitive, suit la loi r\Žflexive car si $a$ \Žtait incomparable \ˆ tous les autres, il faudrait n\Žcessairement postuler $a = a$.

\

4.--{\sc Ordre et classification engendr\Žs par une structure $(M \leqs)$.}

\

Si $a$ et $b$ \Žtaient des \Žl\Žments, distincts ou identiques, d'une paire sym\Žtrique, nous \Žcririons $a = b$. Ainsi, on obtient une structure binaire $(M =)$ qui est une classification. En effet: elle b\Žn\Žficie de la loi transitive, puisque les \Žl\Žments sym\Žtriques \ˆ un tiers sont sym\Žtriques entre eux; le caract\re r\Žflexif est une exigence \Žvidente de la construction; que les paires asym\Žtriques manquent est clair car si nous \Žcrivons  $a = b$ c'est parce que les deux sont vrais
$$a  \leqs b \ \et b\leqs a$$
ce qui nous permet \Žgalement d'\Žcrire $b = a$.

Si $a = b$  et $a\leqs x\leqs b$ sont vraies, alors $a = x = b$ sera \Žgalement vraie.

En effet : la propri\Žt\Ž c) des structures transitives permet d'\Žcrire
$$b\leqs x \leqs a$$
et, d'apr\s la construction de la structure $(M =)$, on peut \Žcrire
$$a = x = b$$

Par cons\Žquent, les classes d'\Žl\Žments \Žgaux sur chaque ligne com\-pl\te
de $(M \leqs)$ sont des intervalles, c'est-\ˆ-dire des ensembles qui, contenant
les \Žl\Žments $a$ et $b$, contiennent tous les interm\Ždiaires, ceux qui v\Žrifient
$$a \leqs x \leqs b$$

Si dans $(M \leqs)$ les \Žl\Žments $a$ et $b$ constituaient une paire asym\Žtrique et que $a \leqs b$ \Žtaient vraie, on \Žcrirait $a <b$.

La structure binaire $(M <)$ ainsi obtenue, \ˆ partir de $(M \leqs)$, est un ordre.

En effet : il est clair que la loi transitive est v\Žrifi\Že, puisque
$$a <b \ \et b < c \ \text{impliquent pr\Žcis\Žment} \ a\leqs b \ \et \ b\leqs c$$
ce qui permet d'\Žcrire $ a \leqs c$. \`A pr\Žsent, on ne peut pas \Žcrire $c \leqs a$  car avec le dernier des ant\Žpr\Žc\Ždents on aurait $b \leqs a$, contraire \ˆ l'hypoth\se que la paire $(a, b)$ est asym\Žtrique; alors nous pouvons \Žcrire $a < c$. Il est clair qu'il n'y a pas d'\Žl\Žments r\Žflexifs car dans $(M \leqs)$ les paires $(a, a$) sont manifestement sym\Žtriques.

Les structures $(M =)$ et $(M <)$ d\Žduites de $(M \leqs)$ sont mutuellement li\Žes selon la loi suivante :
$$a '= a < b \ \text{implique} \  a' < b \ \text{et aussi} \  a < b = b' \ \text{implique} \ a < b'.$$

En effet:  $a '= a$ n\Žcessite $a' \leqs a$, et $a < b$ n\Žcessite $a\leqs b$; des deuxi\me et quatri\me, on d\Žduit $a' \leqs b$. Si l'on avait $b \leqs a'$, jointe  \ˆ $a' \leqs a$, cela donnerait $b \leqs a$, contre $a < b$; donc $a '<b$. La deuxi\me implication se d\Žmontre de mani\re analogue.

De cette propri\Žt\Ž d\Žcoule que, si $a_1$ et $b_1$ \Žtaient deux classes de $(M =)$ et entre deux \Žl\Žments $a$ et $b$ de celles-ci, $a <b$ \Žtait vaie, on aurait
de m\me, pour toute autre paire, $a' < b'$ : cela permet de d\Žfinir l'ordre $(C <)$ entre les classes d\Žfinies par la classification $(M =)$, qui serait isomorphe \ˆ celle subordonn\Že par
$(M <)$ sur l'ensemble $P$ qui ne contiendrait qu'un seul \Žl\Žment de chaque classe. La correspondance qui assigne  \ˆ chaque \Žl\Žment de $M$ la classe \ˆ laquelle il appartient serait un homomorphisme entre $(M <)$ et $(C <)$.

Quand la structure binaire transitive, \Žtant r\Žflexive, manque de paires
sym\Žtriques, dans la structure $(M =)$, d\Žduite de $(M \leqs)$, chaque \Žl\Žment serait unique dans sa classe. La
structure $(M <)$ s'obtiendrait, \ˆ partir de $(M \leqs)$, en changeant le signe $\leqs$ en $<$ dans les paires $a \leqs b$
d'\Žl\Žments diff\Žrents, et en omettant de mettre le signe $<$ entre les paires $a \leqs a$. Dans ce sens, on peut consid\Žrer les ordres comme le cas 221 des structures $(M \leqs)$, selon ce que certains ont coutume de faire.

\

\head{CHAPITRE II. LES ORDRES}

\ 

5.--{\sc Notions pr\Žliminaires.}

\

L'ordre est la structure binaire transitive $(M <)$ d\Žpourvue d'\Žl\Ž\-ments r\Žflexifs donc \Žgalement de paires sym\Žtriques. Le signe $<$, en s'appuyant sur l'intuition, nous le lirons   \guil pr\Žc\de\guir.

Un exemple, assez significatif, est donn\Ž par la figure 3. $M$ est l'ensemble des points du plan contenus dans les lignes de la figure. La relation $a <b$ signifie que l'on peut aller, du premier au second, en suivant toujours  le cours des fl\ches. La paire $b c$ est incomparable.

\

\includegraphics[width=4.8in]{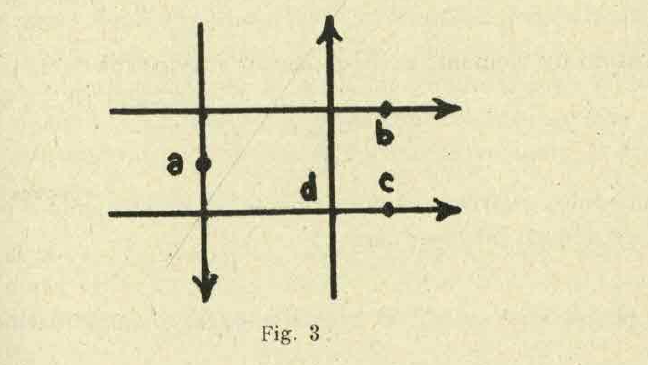}

\

Parmi les structures d'ordre $(M <)$, notons les ordres {\it totaux} dont l'\Žtude a \Žt\Ž entam\Že par Cantor, dans lesquels les paires incomparables manquent. \`A ceux qui ne sont pas totaux, leur inventeur Hausdorff les a appel\Žes {\it partiels}. En parlant d'ordre, sans
pr\Žciser davantage, nous nous r\Žf\Žrons aux uns et aux autres indistinctement.

Il conviendra d'adjoindre, \ˆ la structure $(M <)$, la structure $(M \leqs)$ r\Žsultant de l'\Žcriture $a \leqs b$,  aussi bien si $a <b$, que si $a$ et $b$ \Žtaient identiques.

Dans un ordre, nous appellerons {\it maximum} un \Žl\Žment $a$ lorsque $\ovs a$ est vide; c'est-\ˆ-dire que la relation $a < x$ n'a pas de solutions dans $M$. De la m\me mani\re, un \Žl\Žment $a$ est un {\it minimum} lorsque $\uds a$ est vide; c'est-\ˆ-dire que la relation $x < a$ n'a pas de solution dans $M$. Dans la figure 4, les \Žl\Žments $a$ et $b$ sont tous deux des maximums. 

\

\includegraphics[width=4.8in]{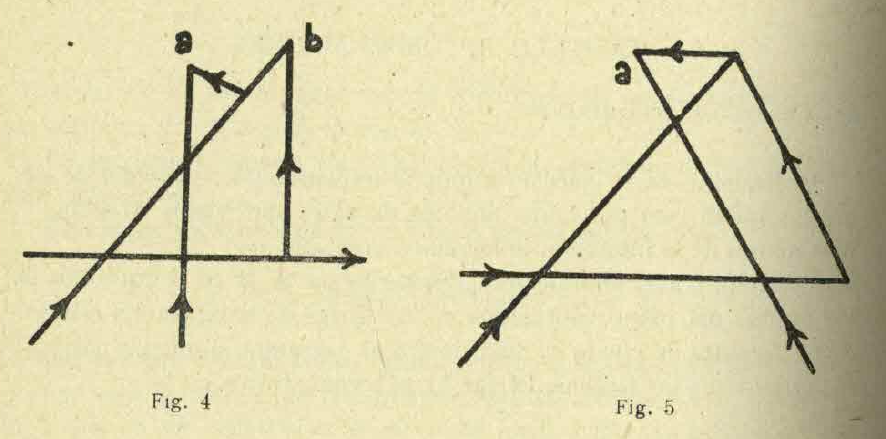}

\noi Les maximums, quand il y en a plusieurs, constituent une transversale.

Lorsqu'un \Žl\Žment $a$ est apr\s tous les autres, c'est-\ˆ-dire
$$\uds a \equiv M-a$$
nous l'appellerons un {\it supremum}. De m\me si l'\Žl\Žment $a$ pr\Žc\de tous les autres, c'est-\ˆ-dire, que
$$\ovs a \equiv M-a$$
nous l'appellerions  {\it infimum}.

Bien s\žr, s'il y a un supremum, c'est le seul maximum, bien qu'il puisse y avoir un maximum unique, sans que cela suffise \ˆ en faire le supremum.

Dans la figure 5, l'\Žl\Žment $a$ est supremum.
\'Etant donn\Ž les ordres $(M <)$ et $(N <)$,  nous appellerons $(M\times N <)$ l'ordre produit des deux, en disposant les \Žl\Žments de l'ensemble $M \times N$, les paires $(m, n)$, lexicographiquement par rapport aux deux.

La con-fusion de structures binaires transitives est toujours possible quel que soit le nombre de structures binaires transitives soumises \ˆ la fusion conjointe; donc, si tout ou partie de celles soumises \ˆ l'op\Žration \Žtaient des ordres, la structure r\Žsultant de leur fusion sera binaire transitive, mais ce ne sera pas toujours un ordre. Par exemple, si les traits de la figure 1 \Žtaient des ordres dans le sens des fl\ches, la structure de con-fusion ne l'est pas en ayant des \Žl\Žments r\Žflexifs. Nous nommerons ordinalement inconfusibles (in-con-fusibles) les  ordres d'un syst\me, lorsque la structure binaire transitive,  r\Žsultat de leur confusion, n'est pas un ordre.

On appelle {\it cycles} d'une structure binaire transitive les syst\mes d'\Žl\Žments qui, pour un certain $a$, v\Žrifient la relation
$$a \ \to \ x \ \to \ a$$
\ˆ condition qu'il existe des solutions autres que $a$ dans le syst\me.

C'est une condition n\Žcessaire et suffisante pour que la structure binaire transitive
r\Žsultant de la con-fusion des ordres soit aussi un ordre, que l'absence de cycles.

En effet : si l'ordre, par d\Žfinition, n'a pas d'\Žl\Žments r\Žflexifs, cela  rend l'existence de cycles impossible;  alors la condition est n\Žcessaire. \`A son tour, si ce n'\Žtait  pas un ordre, il y aurait des \Žl\Žments  r\Žflexifs; supposons que $a$ soit l'un d'entre eux, ainsi $a \ \to \  a$ serait vraie dans la structure r\Žsultant de la fusion conjointe; mais, dans les ordres \Žl\Žmentaires, la relation $a \ \to \  a$ est impossible; alors cela ne peut r\Žsulter qu'en appliquant, apr\s la fusion, la loi
transitive sur l'ensemble fini des relations 
$$a \ \to \ b_1 \ , \ b_1 \ \to \ b_2\ , \ b_2 \ \to b_3\ , \ \dots \ , \ b_n \ \to \ a$$
tir\Žes des ordres  soumis \ˆ la fusion conjointes; alors, on aura un cycle.

On notera que la pr\Žsence de paires sym\Žtriques est \ˆ l'origine de cycles dans la structure; ainsi, leur absence est une condition, n\Žcessaire et suffisante, selon le th\Žor\me pr\Žc\Ždent, \ˆ la confusion des ordres.

En appliquant \ˆ $(M <)$ l'op\Žration que nous avons appel\Že sousjonction, on aura, comme on le voit imm\Ždiatement, un ordre $(M_1 <)$ dans le syst\me des sous-ensembles. Parfois, seules des sous-structures de $(M_1 <)$ seront int\Žressantes, certains sous-ensembles, comme les intervalles, les sections initiales, que nous d\Žfinirons par la suite. 

\newpage

6.--{\sc Majorant et minorant d'un ensemble.}

\

Soit $A$ un sous-ensemble de $M$  qui inclut les solutions de la relation $A <x$, nous l'appellerons le majorant de $A$, et nous le d\Žsignerons par $\ovs A$. De la m\me mani\re, le minorant de $A$ est d\Žfini  que nous d\Žsignerons par $\uds A$.

Il est clair que $\ovs A$ \Žtant l'ensemble des solutions communes au syst\me de relations
$$a < x \ \text{pour} \ a\in A$$
on aura 
\[\ovs A \equiv \prod_{a \in A} \ovs a \ \text{et de m\me} \ \uds A \equiv \prod_{a \in A} \uds a \tag*{[1]}\] 

De la m\me d\Žfinition, suivent imm\Ždiatement les relations
\[\uds A  <  A  < \ovs A\tag*{[2]}\]
\[A \inq B \ \text{implique} \  \ovs B \inq \ovs A \ \et \uds B \inq \uds A\tag*{[3]}\]

En effet :
\[\ovs B \equiv \prod_{b\in B}\ovs b\equiv \prod_{a\in A}\ovs a \ .\prod_{c\in B-A}\ovs c \inq \prod_{a\in A}\ovs a \equiv \ovs A\]

L'ensemble $\ovs A$ est le plus grand parmi ceux qui v\Žrifient  $A < H^i$, donc la r\Žunion de tous.

En effet :
$$A  < H^i  \ \text{implique} \ A < \prod H^i \inq \ovs A$$
de plus, par d\Žfinition,  $A < \ovs A$ ; alors $\ovs A$ est l'un des ensembles $H^i$ donc
$$\ovs A \inq \sum H^i$$
des deux r\Žsulte
\[\ovs A \equiv \sum H^i\tag*{[4]}\]

Sont \Žgalement vraies
\[\ovs A\uds)\ovs) \equiv \ovs A  \et  \uds A\ovs)\uds)  \equiv  A\tag*{[5]}\]

En effet:
\[\ovs A\uds) < \ovs A \ \ \therefore \ \ \ovs A \inq \ovs A\uds)\ovs)  \tag*{[a]}\]
de plus
\[A < \ovs A   \ \text{donne } \ A \inq \ovs A\uds) \ \ \therefore \ \ \ovs A\uds)\ovs) \inq \ovs A\tag*{[b]}\]
De  (a) et (b), tous deux,  r\Žsulte le premier du groupe (5).

Si $B \ . \ \ovs A$ n'est pas vide, nous dirons que $B$ est {\it finalement sup\Žrieur} \ˆ $A$, en \Žcrivant $A <_1 B$.

Cette relation binaire entre les sous-ensembles de $M$ donc manifestement transitive et irr\Žflexive, ordonnera -- partiellement en g\Žn\Žral -- le syst\me $M$ des sous-ensembles de $M$; nous obtenons, ainsi, l'ordre $(M_1 <_1)$.

Pour deux sous-ensembles $A$ et $B$ de $M$, nous \Žcrirons 
$$A =_1  B \ \text{lorsque} \  \ovs A \equiv \ovs B$$

Il y a ainsi une classification $(M_1 =_1)$ pour le syst\me de sous-ensembles de $M$
car il est clair que cette relation binaire est transitive, sym\Žtrique et r\Žflexive.

L'ensemble  r\Žunion de tous ceux d'une classe de $(M_1 =_1)$ appartient \ˆ la m\me classe.

En effet : soient $A_i$, pour $i\in I$, ceux d'une classe; donc, quel que soit $i$, on aura $\ovs{A_i} \equiv P$. Posons
\[\sum_{i\in I} A_i \equiv S \ \ \therefore \ \  A_i \inq S \ \ \therefore \ \ \ovs S \inq \ovs {A_i} \equiv P\tag*{[a]}\]
Comme l'\Žl\Žment g\Žn\Žrique de $P$ d\Žpasse le g\Žn\Žrique $a_i$  de $A_i$,  on aura
\[P \inq \ovs S\tag*{[b]}\]
De (a) et (b), tous  deux,  r\Žsulte
$$\ovs S \equiv P$$
qui stipule que l'ensemble $\displaystyle\sum_{i\in I} A_i$  est dans la m\me classe que ses sommants.

\

7 .--{\sc Fermetures initiale et finale d'un ensemble.}

\

L'ordre $(M <)$ \Žtant donn\Ž,  y {\it fermer initialement} l'ensemble $A$  signifiera faire l'ensemble $\udl A$ qui, avec l'\Žl\Žment $a$, contient tous les \Žl\Žments qui lui sont ant\Žrieurs. La {\it fermeture finale}, dont la d\Žfinition est analogue, sera d\Žsign\Že par $\ovl A$. Comme on l'a d\Žj\ˆ vu au chapitre $1.^\circ$, nous aurons
$$\udl A \equiv A + \sum \uds a \equiv \sum \udl a \ \ \ \ \ \ \ \ovl A \equiv A + \sum \ovs a \equiv \sum \ovl a$$
$a$ \Žtant l'\Žl\Žment g\Žn\Žrique de $A$.

Ces op\Žrations constituent, sur l'ensemble $M$, deux structures topologiques d\Žtermin\Žes par la structure d'ordre $(M <)$.

Les ensembles  ferm\Žs initialement, nous les appellerons {\it sections initiales} de l'ordre $(M <)$.  Les  ferm\Žs finalement, {\it sections finales}. Nous appellerons {\it intervalle} de $(M <)$  les sous-ensembles $S$ de $M$ dans lesquels, \Žtant donn\Žs $a < b$, tous ceux compris entre eux apparaissent \Žgalement dans $S$, c'est-\ˆ-dire ceux qui v\Žrifient
$$a <x < b$$
Si $S$ n'avait pas de paires d'\Žl\Žments comparables, nous l'appellerions \Žgalement un intervalle.

Evidemment  $\udl A \ . \ \ovl A$ sera la fermeture segmentaire de l'ensemble $A$. 
\[\text{De} \ A \inq B \ \text{d\Žcoule} \   \udl A\inq \udl B  \et \ovl A\inq \ovl B\tag*{[7]}\]
En effet : 
$$\udl B \equiv \sum_{b\in B} \udl b \equiv \sum_{a\in A} \udl a + \sum_{c \in B - A} \udl c \ \ \therefore \ \ \udl A \equiv \sum_{a\in A} \udl a \ \inq \ \sum_{b\in B} \udl b \equiv \udl B$$
qui est le premier du groupe (7).

Une cons\Žquence imm\Ždiate de la d\Žfinition de la section initiale est  que la r\Žunion, et aussi l'intersection, d'un syst\me de sections initiales, est encore une  section initiale.

En effet : soient $A_i$ pour $i \in I$  les sections initiales du syst\me consid\Žr\Ž. $\displaystyle\sum_{i\in I} A_i$ contient, avec les \Žl\Žments $a \in A_i$, tous ceux de l'ensemble $\uds a$;  c'est donc une section initiale. De m\me, si $a \in\displaystyle \prod_{i\in I} A_i$, c'est qu'il appara\"t dans tout $A_i$, quel que soit $i\in I$; alors il appara\"t \Žgalement dans tous les $A_i \uds a$, et donc $\displaystyle\prod_{i\in I}A_i$ est aussi section initiale.

La fermeture initiale $\udl A$ d'un ensemble $A$, dans $(M <)$, est l'intersection des sections initiales contenant $A$.

En effet : soit S la section initiale g\Žn\Žrique dans laquelle $A$ est contenu
$$\text{de} \ A \inq S \ \text{vient} \  \udl A \inq \udl S = S$$
Comme $\udl A$ est \Žgalement l'une des sections initiales contenant $A$, il en r\Žsultera ce qui pr\Žc\de, et que $\udl A$ est le minimum de celles qui contiennent $A$.

On dit qu'un ordre $(M <)$ est {\it ramifi\Ž}, (Kurepa) quand $(\udl m <)$ est un ordre total, pour tout $m \in M$. La figure 6 donne un exemple de ces ordres tr\s int\Žressants.

\

\

\includegraphics[width=4.8in]{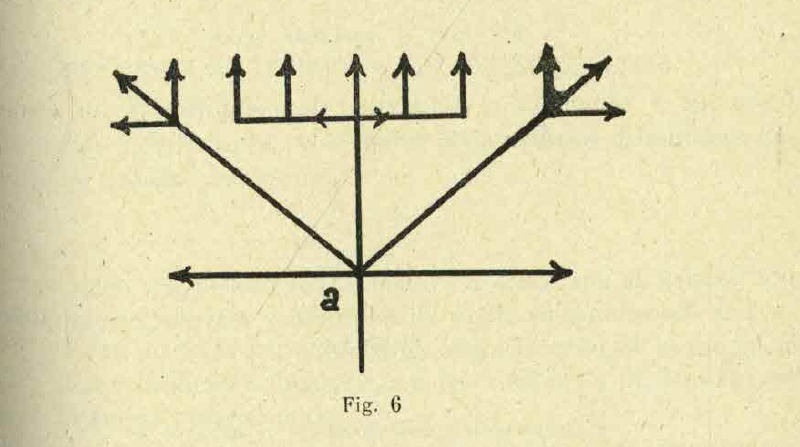}

\

\

Nous disons {\it cofinaux} dans $(M <)$ les ensembles $A$ et $B$ et, plus pr\Žcis\Žment les ordres induits, lorsque leurs fermetures initiales co\•n\-cident;
c'est \ˆ dire
$$(A <) \et \ (B <) \ \text{sont cofinaux} \  \si  \udl A \equiv \udl  B$$
Corr\Žlativement,  sont d\Žfinis les {\it co\•nitiaux}.

Chacune de ces relations d\Žfinit, sur $M$, une classification; c'est-\ˆ-dire, en nous r\Žf\Žrant \ˆ ce qui a \Žt\Ž express\Žment nomm\Ž, nous \Žcrirons
$$A \equiv_2 B \ \text{lorsque} \ \udl A \equiv \udl B$$

C'est une condition n\Žcessaire et suffisante pour la confinalit\Ž de $A$ et $B$  que, quelle que soient
 ses \Žl\Žments g\Žn\Žriques, $a$ et $b$ respectivement, la relation 
 $$a \leqs x \ \text{ait une solution dans} \ B, \et b\leqs x\ \text{ait une solution dans} \ A.$$

 En effet : supposons qu'ils soient cofinaux. Ayant un $a\in \udl A \equiv \udl B$, soit il est de $B$, soit il est d\Žpass\Ž, donc appartient \ˆ $\udl B$, par certains de $B$. Identiquement pour le second. Ainsi, la condition est n\Žcessaire. Supposons maintenant que ces relations aient les dites solutions; par la
premi\re, $\udl A \inq \udl B$; par la seconde, $\udl B \inq \udl A$. Des deux vient $\udl A \equiv \udl B$, qui est la condition d\Žterminante de la cofinalit\Ž.

Nous dirons qu'un ensemble $B$ {\it enveloppe sup\Žrieurement} un autre $A$, dans l'ordre $(M <)$, lorsque,
$$\ovl a . B \not\equiv {\rm O}  \pourtout a \in A$$
c'est-\ˆ-dire que $B$ p\Žn\tre dans la fermeture finale de l'un quelconque des \Žl\Žments de l'ensemble $A$.
Nous \Žcrirons, alors
$$A \leqs_3 B$$
parce que c'est une relation \Žvidemment transitive et r\Žflexive.

Les op\Žrations de fermeture initiale (finale) et de majoration (minoration) sont li\Žes par les
propositions suivantes:
\[\ovs A \ \text{est une section finale de} \ (M <), \text{c'est-\ˆ-dire} \  \ovs A{^-} \equiv  \ovs A \tag*{[8]}\]

En effet : si $p\in \ovs A$, c'est parce que 
$$A < p < \ovs p \ \ \therefore \ \ \ovs p \ \inq \ovs A$$
c'est-\ˆ-dire que $\ovs A$ contient, avec $p$, tous ceux qui le suivent.

C'est la m\me chose que de majorer $A$ ou sa fermeture initiale; c'est-\ˆ-dire
\[\ovs A \equiv \udl A)^* \tag*{[9]}\]

En effet :
\[\text{De} \  A \ \inq \ \udl A\ \ \text{vient} \ \ \udl A)^* \ \inq \ \ovs A\tag*{[a]}\]
Quel que soit l'\Žl\Žment $q$ de $\udl A$, il sera surpass\Ž, ou co\•ncidera, avec certains $a$ de $A$. En tant qu'\Žl\Žment g\Žn\Žrique $x$ de $\ovs A$ v\Žrifie $a < x$, quel que soit $a \in A$; alors, dans les deux cas, $q < x$; donc $x$ appartiendra \ˆ l'ensemble $\udl A)^*$ , et ainsi
\[\ovs A \ \inq \ \udl A)* \tag*{[b]}\]
De (a) et (b), tous deux, resulte (9).
\[ \udl A \equiv \udl B \ \text{implique} \ \overset {*}A \equiv \ovs B\tag*{[10]}\]

En effet :
$$\ovs A \equiv  \udl A) ^* \equiv \udl B)^* \equiv \ovs B$$

\

8.--{\sc Extension syst\Žmatique des ordres.}

\

Entre deux ordres $(M <_1)$ et $(P < _2)$ nous disons que celui-ci est une {\it extension} de celui-l\ˆ  et nous \Žcrivons
$$(M <_1) \leqs' (P < _2)$$
lorsque :

1.$^{\circ})$ $P$ est un sur-ensemble de $M$.

2.$^{\circ}$) Les relations d'ordres, pour les \Žl\Žments de $M$, sont les m\mes dans les deux ordres.

Le signe $\leqs'$, pour aider l'intuition, nous le lirons \guil immerg\Ž dans\guir 

 De l'ordre $(M <_1)$ nous dirons que c'est un {\it sous-ordre} de $(P < _2)$; de celui-ci que c'est un {\it sur-ordre} de $(M <_1)$.

 Nous devons noter que, \Žtant donn\Ž $(P < _2)$, en donnant simplement le sous-ensemble $M$ de  $P$,  on d\Žtermine $(M <_1)$; par cons\Žquent, dans ce cas, on peut utiliser le m\me symbole
du d\Žpart, dans le construit. Au contraire, se donner $(M <_1)$ et le sur-ensemble strict $P$ de $M$,  ne suffisent pas pour d\Žterminer $(P < _2)$, extension de celle-l\ˆ, mais, en g\Žn\Žral, la relation
\[(M <_1) \leqs' (P <_x)\tag*{[11]}\]
poss\de plusieurs solutions, en l'inconnue $x$. 
Pr\Žcis\Žment le probl\me de la construction syst\Žmatique des ordres inclut celui de l'obtention de toutes les solutions de la relation (11).

Une notion importante est celle d'{\it extension limite}. Supposons que le syst\me des ordres
$$(M_i  <_i ) \ \pour i\in I$$
s'\Žtende en suivant la variable $i$, d\Žterminant le syst\me, dans un sens croissant, l'ordre total, ouvert sup\Žrieurement, $(I <)$; d\Žsignant $M$ l'ensemble limite, nous appellerons {\it extension limite}, et nous \Žcrirons
$$(M <) \equiv \underset{i\in (I <)} {\overset{\longrightarrow}\lim} (M_i <_i)$$
l'ordre $(M <)$, lorsque c'est une extension de tout le syst\me.

Une structure que nous soup\c connons d'\tre importante dans l'\Žtude des ordres est celle constitu\Že, \ˆ travers la relation transitive et r\Žflexive $\leqs'$ \guil immerg\Že dans\guir, entre tous les ordres
$$(S_i <_{i_j}) \ i\in I \ i_j \in J_i$$
du syst\me $S_i$ de tous les sous-ensembles d'un $M$ donn\Ž, homomorphe \ˆ la structure d'inclusion $(M_1 \inc)$ d\Žfinie pour tous les sous-ensembles de $M$.

\'Etant donn\Že une ligne compl\te $(L_1 \inc)$ de $(M_1 \inc)$, l'ensemble des ordres possibles, sur la totalit\Ž des sous-ensembles  qui la constituent, articul\Ž par la relation stricte d'\guil \ immersion \guir \ est manifestement un ordre ramifi\Ž, puisque le sous-ordre de l'une donn\Že est d\Žtermin\Ž de mani\re unique.

En vertu du th\Žor\me du bon ordre -- ressource puissante, sans son aide peu de choses seraient accomplies dans l'\Žtude g\Žn\Žrale des structures -- l'ensemble $(P - M)$ peut \tre dispos\Ž, de plusieurs mani\res,  en une suite bien ordonn\Že
$$p_0, \ p_1, \ p_2, \ \dots \  (p_\alpha$$
et, \Žtant donn\Ž l'ordre $(P <_2)$, on peut consid\Žrer les sous-ordres d\Žrermin\Žs par les ensembles
$$M , M+p_0, M+p_0 + p_1, M+p_0 + p_1 +p_2 \dots$$
donc on peut former l'ordre $(P <_2)$  partant de $(M <_1$) en ajoutant les nouveaux \Žl\Žments un \ˆ un, de mani\re convenable, et en passant \ˆ la limite lorsqu'une suite ouverte d'entre eux a \Žt\Ž ajout\Že.

L'extension limite n'impliquant pas de nouvelles relations ordinales, puisque seules celles obtenues aux \Žtapes pr\Žc\Ždentes sont adjointes, il suffit d'examiner les diff\Žrentes mani\res possibles d'\Žtendre un ordre donn\Ž, par l'ajonction d'un seul \Žl\Žment.

Adjoignant un \Žl\Žment $p$, non inclus dans l'ensemble $M$, et formant sur $M + p$  un ordre $(M + p  <)$, qui laisse les relations ordinales qui liaient d\Žj\ˆ les \Žl\Žments de $M$ invariants, dans l'ordre $(M <_1 )$, on distingue, dans $M$, trois classes d'\Žl\Žments: ceux d'un ensemble $A$, qui en $(M + p <)$ constituent l'ensemble $\uds p$, ceux d'un autre ensemble $B$, qui en $(M + p <)$ constituent $\ovs p$, et le reste, ceux de $M - (A + B)$, incomparables \ˆ $p$,
dans l'ordre $(M + p <)$. On a montr\Ž que $\uds p$ est une section initiale de $(M <_1)$, et $\ovs p$ une autre finale,  les deux v\Žrifiant
$$\uds p \ <_1 \ \ovs p$$
Ainsi, ce proc\Žd\Ž, conduit \ˆ la notion importante suivante :

{\it Trou} d'un ordre $(M <_1)$, nous appellerons la paire $(A, B)$, constitu\Že d'une
section initiale $A$ et une autre finale $B$, qui v\Žrifient $A <_1 B$. (*)

\

\noi  --------------

(*)  Cette notion, g\Žn\Žralisant celle  du m\me nom, que, seule pour les ordres totaux, nous avons introduite dans notre th\se : {\sl Rev. Mat. Hisp.-Am.} (IV) {\it 3 }(1943). D\Žsormais, elle vaut \Žgalement pour les ordres partiels.

\noi ---------------

\

 On dit que le nouvel \Žl\Žment $p$ occupe le trou $(A, B)$, dans l'extension $(M+ p <)$ de $(M <_1)$ quand les relations d'ordre de $p$ avec les \Žl\Žments
de $M$ sont d\Žfinis par la double relation 
$$A < p < B$$
 ceux de $M - (A + B)$ \Žtant incomparables \ˆ $p$; donc 
 $$\uds A  \equiv p \ \et B \equiv \ovs p$$
dans l'ordre construit

On a d\Žj\ˆ  dit que, dans toute extension d'un ordre, avec un nouvel \Žl\Žment, celui-ci 
occupe un trou. Il nous reste \ˆ  montrer que, quel que soit le trou dans l'ordre $(M <_1)$, il existe une extension $(M + p <)$ dans laquelle le nouvel \Žl\Žment occupe  ce trou.

En effet : \Žtant donn\Ž le trou $(A, B)$ de $(M <_1)$, adjoignons \ˆ ses relations binaires  celles qui r\Žsultent en posant
$$A  < p  < B$$
La structure binaire ainsi obtenue est manifestement transitive lorsque $p$ n'intervient pas dans les paires du transit, ou lorsqu'il ne s'agit pas d'un moyen terme; quand c'en est aussi, car $A  <_1 B$. Qu'il n'y ait pas d'\Žl\Žments r\Žflexifs est \Žvident, car ceux de $M$ ne le sont pas, par hypoth\se, ni $p$ par construction. Ainsi $(M + p <)$ est un ordre, et en lui $p$ a occup\Ž le trou $(A, B)$.

Par cons\Žquent, {\it les trous, de l'ordre $(M <_1)$, indiquent toutes les extensions possibles, par adjonction d'un nouvel \Žl\Žment.}

Afin de construire syst\Žmatiquement toutes les structures d'ordre possibles sur un ensemble $M$, on formera une suite bien ordonn\Že  avec ses \Žl\Žments, et
en commen\c cant par le premier, les \Žl\Žments sont adjoints un \ˆ un, en pla\c cant l'adjonction, de toutes les mani\res possibles, dans l'ordre d\Žj\ˆ r\Žalis\Ž, c'est-\ˆ-dire en occupant successivement tous les trous de l'ordre mentionn\Ž. Nous avons d\Žj\ˆ utilis\Ž ce proc\Žd\Ž, de formation syst\Žmatique, pour les totaux dans notre th\se mentionn\Že. Le lecteur averti notera l'analogie de ce proc\Žd\Ž -- mutatis mutandis  -- avec celui suivi par Steinitz, dans son \guil \ Algebraische Theorie der K\šrper  \ \guir, pour la construction exhaustive des corps.

Si $r$ et $s$ appartiennent \ˆ $M - (A + B)$, lui  appartiennent aussi les solutions de la  double relation
$$r <_1 x <_1 s$$
car $x <_1 A$ impliquerait l'inclusion de $r$ dans $A$; de m\me $x < _1 B$ impliquerait celle de $s$ dans $B$. Par cons\Žquent, l'ensemble $M - (A + B)$ est un intervalle de l'ordre $M <_i)$ : nous l'appellerons {\it l'intervalle neutre} du  trou $(A, B)$.

Nous devons noter que, bien que $(M <_1)$ ne soit pas connexe, on pourra \Žtablir, au moyen de l'\Žl\Žment joint $p$, la connexion entre les diff\Žrentes parties qui le composent. 

\

9.--{\sc Les diff\Žrents types de trous.}

\

Nous appelons  trou {\it ext\Žrieur} celui, $({\rm O,O})$, dont l'intervalle neutre est tout $M$.

Nous appelons  trous {\it couverts} ceux, $({\rm O}, B)$, dont la section initiale est vide. De mani\re analogue, trous {\it appuy\Žs } ceux, $(A,{\rm O})$, dont la section  finale est vide.

Trous {\it internes}  appelons-nous les autres, pour lesquels ne sont  pas vides,  ni  leur section initiale, ni leur section finale.

Il est facile de voir qu'il y a des trous $(A, B)$, pour lesquels il n'y a pas de solution dans $M$ \ˆ la
double relation
$$A <_1 x  <_1 B$$

Par exemple, dans la figure 3; le trou $(\udl d, \ovl b + \ovl c)$.

Plus int\Žressants sont les trous $(A, B)$, pour lesquels $x <_1 B$ n'a pas de solution dans $M - A$, qui sont les trous $(\uds B, B)$. 
Leur sont analogues les trous $(A, \ovs A)$. Il y a encore  l'interf\Žrence des deux circonstances, c'est-\ˆ-dire des trous $(A, B)$ pour lesquelles
$$A \equiv  \uds B \ \et B \equiv \ovs A$$
que nous appelons  {\it \Žtroits}. Selon les \Žquations (5) du num\Žro 6, 
ces trous  sont ceux de la forme
$$\left[\uds B, \uds B)^*\right] \ \et  \left[\ovs A)_*, \ovs A\right]$$
La raison du nom est que son intervalle neutre est le plus petit possible, dans le sens que l'on ne peut pas passer de ses \Žl\Žments, adjoints \ˆ  la section initiale ou finale, de mani\re que forment toujours un trou les nouveaux ensembles $A'$, $B'$.

Un trou \Žtroit $(A, B)$ est {\it disjonctif} lorsque son intervalle neutre est vide, consommant ainsi, entre $A$ et $B$ tous les \Žl\Žments de $M$. Dans la figure 6, le trou $(\udl a,\ovs a)$ est disjonctif. Dans les ordres correspondants aux figures 3, 4, 5, il n'y a pas de trous disjonctifs, bien  qu'il y ait des trous \Žtroits.

Dans les ordre totaux, lorsque   l'extension pr\Žvue est \Žgalement totale, les seuls trous qui pr\Žsentent un int\Žr\t sont les \Žtroits, qui sont \Žgalement tous disjonctifs.

Les trous disjonctifs, \Žventuels dans un ordre, ont la propri\Žt\Ž suivante : 

D\Žsignant par $A _i$, pour $i \in I$,  la partie  initiale g\Žn\Žrique d'un syst\me de trous
disjonctifs, de l'ordre 
$$(M <_i),   \ \sum_{i\in I} A_i \ \et \prod_{i\in I}A_i$$
sont aussi des sections initiales de trous disjonctifs.

En effet : comme le montre le no. 7, l'un et l'autre, sont des sections initiales de l'ordre $(M <_1)$. Soit alors $r$ un \Žl\Žment non inclus dans $\displaystyle\sum_{i\in I} A_i$, 
de sorte que, n'apparaissant dans aucun des $A_i$, il sera dans chaque $B_i$ partie finale du trou  disjonctif $(A _i, B_i)$  alors
$$A_i <_1 r  \pour i\in I$$ avec lequel $r$ appartiendra \ˆ l'ensemble $\left(\displaystyle\sum_{i\in I} A_i\right) ^*$ et, ainsi, de cette mani\re, car $r$  est l'\Žl\Žment g\Žn\Žrique de $M- \displaystyle\sum_{i\in I} A_i$, on aura
$$\sum_{i\in I} A_i < _1  M- \sum_{i\in I} A_i$$
ce qui confirme la premi\re partie de notre proposition. La seconde se d\Žmontre facilement  partant de la relation
$$\prod_{i\in I}A_i \equiv M - \sum_{i\in I} B_i$$
et appliquant la partie d\Žmontr\Že \ˆ l'ordre inverse.

\

10.--{\sc Ordre naturel des trous.}

\

\'Etant donn\Žs $(A, B)$  et $(A', B')$  deux trous dans l'ordre $(M <_1)$,
{\it occuper simultan\Žment} ces trous, respectivement avec les \Žl\Žments $p$ et $q$  c'est d\Žfinir une extension $(M + (p,q) <)$ mettant, comme nouvelles relations fondamentales, en plus de celles de $(M <_1)$ celles implicites dans
$$A < p < B \ \et  A'< q < B'$$
et appliquant la loi transitive, jusqu'\ˆ ce que soit ferm\Ž le syst\me de relations d'ordre que nous avons ainsi.

De mani\re analogue, l'occupation simultan\Že de tout syst\me de trous dans un ordre est d\Žfinie.

Nous d\Žfinirons un ordre naturel, pour les trous de $(M <_1)$, en posant $$(A, B) <_1 ( A', B') \ \text{lorsque} \ BA' \not\equiv 0$$

Ce crit\re ordinal est justifi\Ž, parce que, occupant simultan\Žment tous les trous de $(M <_1)$, chacun avec un nouvel \Žl\Žment, la relation d'ordre entre eux sera pr\Žcis\Žment celle attribu\Že aux trous qu'ils occupent; de plus les relations, entre les nouveaux \Žl\Žments, doivent se d\Žduire \ˆ l'aide de la loi transitive; alors, si les \Žl\Žments $p$ et $p'$, ajout\Žs dans les trous \Žcrits ci-dessus, v\Žrifient $p < p'$, c'est qu'il y avait un \Žl\Žment interm\Ždiaire, qui serait n\Žcessairement dans $B.A'$.

Cette relation binaire, entre les trous de $(M <_1)$, est manifestement irr\Žflexives car
$$A \ . \ B \equiv 0$$

D\Žmontrons qu'elle est transitive
$$(A, B) <_1 (C, D) \ \ \text{implique} \ B \ .\ C   \not\equiv 0$$ 
$$(C, D) <_1 (E, F)  \ \text{implique} \ D \ .\ E   \not\equiv 0$$
soit
$$p\in B.C \et q\in D. E  \ \therefore \ p <_1 q \ \therefore \ p\in E$$
par suite
$$B. E \not\equiv  0  \ \therefore \  (A, B) <_1 (E, F)$$

Ainsi, la relation \Žtablie entre les trous est un ordre,  deux trous \Žtant incomparables lorsque
$$B \ . \ A'\equiv  0 \equiv B'\ . \ A$$
et comparable lorsqu'une seule de ces deux est v\Žrifi\Ž.

Si un trou $(A, B)$ pr\Žc\de un autre $(A', B')$, la section initiale de celui-l\ˆ est incluse dans celui-ci.

En effet : par hypoth\se $B, A '\not\equiv {\rm O}$

\noi soit
$$p\in B \ . \ A' \ \therefore \ p\in B  \ \therefore \ A <_1 p \ \therefore \ A\inq \uds p$$
\Žgalement
$$p\in A' \ \therefore \ \uds p \inq A'$$
donc
$$A \inq \uds p  \inq A' \ \therefore \  A \inq A'$$

La condition $A \inq A'$ ne suffit pas pour que deux trous soient comparables car
nous d\Žmontrerons :

Si
$$A \inq A' \ \et B \inq B'$$
les trous $(A, B)$ et $(A', B')$ sont incomparables.

En effet :
$$p\in B \ \text{donne} \ p \in B' \ \therefore \ p\notin A' \ \therefore \ B A' \equiv 0$$
$$q\in A \ \text{donne} \ q \in A' \ \therefore \ q\notin B' \ \therefore \ A B' \equiv 0$$

\

Si les trous comparables sont tous deux disjonctifs et distincts, il est clair
que
$$A . B ' \equiv 0 \ \text{implique} \  A  \inc A'  \ \therefore \  B .  A' \not\equiv  0$$
alors deux trous disjonctifs sont toujours comparables, et la section initiale du pr\Žc\Ždent inclut l'initiale du suivant; donc l'ordre subordonn\Ž  \ˆ l'ordre naturel des trous, pour les trous disjonctifs, est total.

Cet ordre total, des trous disjonctifs de $(M <_1)$, a des sauts.

En effet : d\Žsignons par $(A <_i, B_i)$ les trous disjonctifs, pour lesquels un \Žl\Žment donn\Ž $m$, de $M$, appartient \ˆ leur section finale; de m\me $(\alpha_j, \beta_j)$ d\Žsigne les trous disjonctifs aux sections initiales desquels appartient $m$. Comme on le montre, \ˆ la fin du no. 9,  sont aussi des trous disjonctifs
$$\left(\sum_{i\in I} A_i, \prod_{i\in I}B_i\right) \ \et \left(\prod_{j\in J} \alpha_j , \sum_{j\in J}\beta_j\right)$$
et, par la d\Žfinition, donn\Že ci-dessus, on  aura
$$(A_i, B_i) \leqs_1 \left(\sum A_i, \prod B_i\right) <_1\left(\sum \alpha_j , \prod \beta_j\right)\leqs_1(\alpha_J, \beta_j)$$
les trous disjonctifs \Žtant \Žpuis\Žs avec les deux classes, les trous form\Žs sont contigus.

Soient $(A'_r, B'_r)$ et $(\alpha'_s, \beta'_s)$ les syst\mes de trous disjonctifs, form\Žs de mani\re analogue, pour un autre \Žl\Žment $m'\in M$ : si
on avait, car  tous sont comparables, \Žtant des trous disjonctifs, 
$$\left(\sum A_i, \prod B_i\right)<_1\left(\prod \alpha_j , \sum \beta_j\right)<_1\left(\sum A'_r, \prod B'_r\right)<_1$$
$$<_1\left(\prod \alpha'_s , \sum \beta'_s\right)$$
on aurait
$$\prod B_i . \prod \alpha_j <_1\sum \beta_j . \sum A'_r <_1 \prod B'_r . \prod \alpha'_s \ \therefore$$
$$\therefore \prod B_i .\prod \alpha_j <_1 \prod B'_r . \prod \alpha'_s$$
Donc, les \Žl\Žments de $M$ sont r\Žpartis dans les ensembles analogues \ˆ  $\prod B_i . \prod \alpha_j$ , deux \ˆ deux comparables, dans l'ordre $(M_1 <_1)$, correspondant au syst\me de sous-ensembles de $M$. Comme il est \Žvident que l'ordre
$$\left[\left(\prod B_i . \prod \alpha_j\right)\right]$$
n'a pas de trous disjonctifs, nous aurons :

Tout ordre $(M <_1)$ est d\Žcompos\Ž, de mani\re unique, en un syst\me totalement ordonn\Ž, au moyen de
la relation $<_1$, des ensembles, et l'ordre induit par $(M <_1)$, sur chacun d'eux, n'a pas de tous disjonctifs.

L'ordre, tir\Ž de la figure 7, donne un exemple d'ordre avec les trous disjonctifs suivants:
$$(\uds a, \bar a) \ (\udl a, \ovs a) \ (\uds b, \bar b)\ (\udl b, \ovs b) \ (\uds c, \bar c) \ (\udl c, \ovs c) \ (\uds d, \bar d) \ (\udl d, \ovs d) \ (\uds e, \bar e) \ (\udl e, \ovs e)$$ 
et sa d\Žcomposition, selon la derni\re proposition, donne les ordres sur les
ensembles, que nous \Žcrivons pr\Žcis\Žment dans l'ordre o\ ils apparaissent
$$\uds a, a, \ovs a \ \uds b, b, \ovs b \ \uds c, c, \ovs c \ \uds d, d, \ovs d \ \uds e, e, \ovs e$$ 

\

\

\includegraphics[width=4.8in]{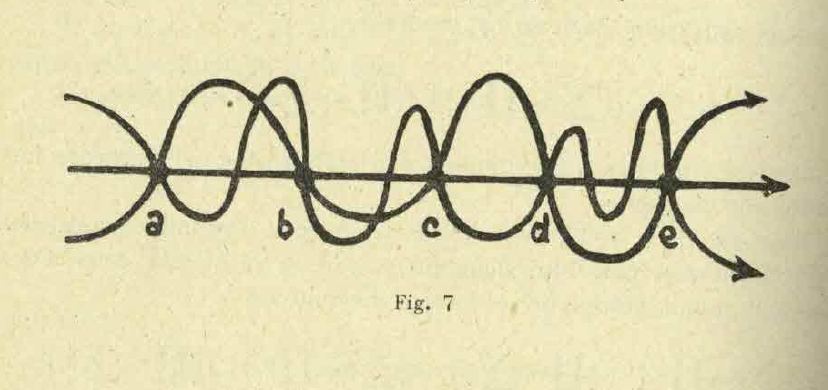}

\

\

11.--{\sc Matrice d'ordre pour un nombre cardinal.}

\

Suivant une d\Žnomination tr\s suggestive de Denjoy (*), nous appellerons ordre {\it matrice}, pour le nombre cardinal $K$, un ordre $(U <)$, lorsque quel que soit l'ordre $(M <')$, v\Žrifiant $|M| = K$, parmi les sous-ordres de $(U <)$, il y en a un 
isomorphe \ˆ $(M <')$.

\

\noi -----------------

(*) {\it L'\Žnum\Žration transfinie} (Livre 1.$^\circ$ - 1946 - page 146). En tant qu'ordres universels, nous les envisageons dans notre th\se, mentionn\Že plus haut, pour les ordres totaux. Maintenant, il faut le noter, nous nous r\Žf\Žrons aux ordres quelconques.

\noi------------------

\ 

Il est facile de d\Žfinir des matrices d'ordre, \ˆ laide du proc\Žd\Ž suivant :

\'Etant donn\Ž un ordre quelconque $(Z <)$, nous consid\Žrerons l'op\Ž\-ration d'extension $\f$ 
consistant \ˆ remplir, simultan\Žment, tous les trous dans $(Z <)$, mais avec un seul objet chacun : on obtiendra ainsi l'ordre $(Z\f <)$. En r\Žp\Žtant successivement, on obtiendra une suite bien ordonn\Že d'extensions, sur les ensembles
$$Z\f_\alpha \equiv (Z\f_{\alpha - 1})\f \pour \alpha \ \text{ordinal non limite}$$
$$Z\f_\alpha \equiv \underset{i < \alpha}{\rm Lim} \ Z\f_i \pour \alpha \ \text{ordinal limite}$$
sous-entendu, pour uniformiser la notation,
$$Z \equiv Z\f_0$$

Consid\Žrant ordonn\Ž l'ensemble d'un objet, nous parlerons de l'ordre $(I <)$. En lui appliquant l'op\Žration d'extension $\f$, les ordres $(I\f_i  <)$ seront form\Žs qui, \ˆ mesure que l'ordinal $i$ cro\"t, sont \Žtendus.

D\Žsignant par $\gamma$  le plus petit ordinal qui v\Žrifie 
$$|\gamma| = |M|$$
je dis que {\it l'ordre $(I\f_{\gamma} <)$ est une matrice des ordres de $M$}.

En effet : selon le th\Žor\me du bon ordre, il existera pour $M$ un
bon ordre $(M <_{b.0.})$ de type ordinal $\gamma$. Soit $m$ le premier \Žl\Žment de celui-ci. En passant de $(m <)$ \ˆ $(M <)$, certains trous sont occup\Žs; en prenant dans $(M <_{b.o.})$ les plus bas qui les occupent, nous obtiendrons $(m\psi_1 <)$. En occupant les autres  trous de $(m <) $ chacun avec un \Žl\Žment n'appartenant pas \ˆ $M$, nous obtiendrons  $(m\f_1 <)$.

Sont \Žvidemment vraies
$$(m\psi_1 <) \ \text{est immerg\Ž dans} \ (M <)$$
$$(m\psi_1<) \ \text{est immerg\Ž dans}\ (m\f_1<)$$
En passant de $(m, <)$ \ˆ $(M <)$ certains trous $(A, B)$ de celui-l\ˆ seront occup\Žs; en prenant de $M$ les plus bas dans $(M <_{b. o.})$ qui les occupent, on aura
$$(m\psi_2 < ) \ \text{est immerg\Ž dans}\  (M <)$$
Fermetures initiale et finale, par rapport \ˆ $(m\f_1 <)$, $A$ et $B$,
on obtient le trou  $(\udl A, \ovl B)$ qui,  \ˆ l'\Žtape $(m\f_2, <)$ est occup\Ž par le m\me \Žl\Žment de $M$ qui occupait auparavant le trou $(A, B)$  : on obtiendra, en remplissant les trous restants chacun avec un \Žl\Žment ext\Žrieur \ˆ $M$
$$(m\psi_2 < ) \ \text{est immerg\Ž dans}\  (m\f_2 < )$$
On r\Žit\re ce proc\Žd\Ž, en le combinant avec l'extension limite. En appliquant l'induction [la r\Žcurrence], on d\Žmontre facilement 
$$(m\psi_\alpha <) \ \text{est immerg\Ž dans}\  (m\f_\alpha < )$$
le premier \Žtant \Žvidemment immerg\Ž dans $(M <)$. Comme l'ensemble $M$ s'\Žpuise avec $m\psi_\sigma$
pour $\sigma \leqs \gamma$, il en r\Žsulte
$$(M <) \equiv (m\psi_\sigma <) \ \text{est immerg\Ž dans}\ (m\f_\sigma <) \ \text{est immerg\Ž dans}\ (m\f_\gamma <)$$
comme nous voulions le d\Žmontrer (*)

\noi------------------

(*) Nous soup\c connont que le cardinal de $m\f_\gamma$ est $2^{|M|}$

\noi -----------------

\

\

12.--{\sc Lignes compl\tes}

\

Pour les structures binaires ordinales, $(M <)$, les lignes, aussi bien compl\tes qu'incompl\tes, semblent d'un int\Žr\t exceptionnel. Le sous-ensemble $L$, avec l'ordre induit $(L <)$, sera une {\it ligne}, si son ordre est total. La ligne sera compl\te, lorsque, quel que soit l'\Žl\Žment $q$ de $(M - L)$ ajout\Ž, $(L + q <)$ est d\Žj\ˆ un ordre partiel.

La notion de ligne peut \tre relativis\Že car en d\Žsignant par $R$ un sous-ensemble de $M$, l'ensemble $LR$ sera aussi une ligne de $(R <)$. La ligne compl\te non; ainsi, sur la figure 8, o\ $L$ \Žtant la ligne, \ˆ l'extr\Žmit\Ž terminale de laquelle le $L$ appara\"t, et
$$R \equiv \udl c + \udl d$$

\

\includegraphics[width=4.8in]{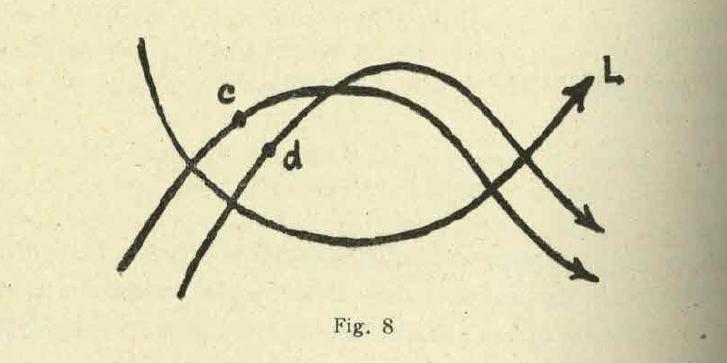}

\

\noi $LR$ est une ligne incompl\te de $(R <)$.

 $(S <)$ \Žtant  une ligne de $(M <)$, il y a des lignes compl\tes qui contiennent $S$.

Ce th\Žor\me a \Žt\Ž d\Žmontr\Ž par Hausdorff. On compl\te $(S <)$ en occupant ses trous disjonctifs $(I,F)$ avec l'une des lignes compl\tes, soit $(L_{if} <)$, de
$(\ovs I\uds F <)$; cela donne toutes les mani\res possibles de compl\Žter $(S <)$.

Parmi les lignes, il convient de signaler les {\it demi-rayons initiaux}, les lignes $(I <)$ telles que, les \Žl\Žments
comparables, \ˆ tous ceux de $I$, appartiennent \ˆ $\ovs I$; donc, $(I <)$  peut \tre compl\Žt\Ž seulement par son trou final $(I,{\rm O})$; ainsi ils pourraient \tre d\Žfinis comme des lignes compl\tes de la structure ordinale $(M - \ovs I <)$. Corr\Žlativement, on d\Žfinit le {\it demi-rayon final}.

\'Etant donn\Že la ligne compl\te $(L <)$  de $(M <)$, et $(I, F)$ \Žtant  l'un des trous disjonctifs de celui-l\ˆ, on v\Žrifie \Žvidemment que $I$ est un demi-rayon initial, $F$ un demi-rayon final de $(M <)$.

De m\me, l'ensemble $\ovs I .\uds F$ est un ensemble vide. 

En effet : si $m$ en \Žtait un \Žl\Žment, on aurait
$$I  <m < F$$
et la ligne, contrairement \ˆ l'hypoth\se, ne serait pas compl\te.

Les demi-rayons $I, F$ engendrent  les deux suites d'ensembles
$$\text{le} \ I  \ \udl I, \ovs I, \overset {*}I\uds);  \ \text{le} \  F \  \ovl F, \uds F, \uds F\ovs)$$
qui v\Žrifient les inclusions suivantes : 
$$I \inq  \udl I\inq \ovs I\uds) \inq \uds F  \et  \  F\inq \ovl F \inq  \uds F\ovs)  \inq \ovs I$$

En effet, et se r\Žf\Žrant uniquement au premier groupe, puisque le second est \Žtabli corr\Žlativement, en utilisant l'ordre inverse : le premier
$$I \inq \udl I$$
est une cons\Žquence \Žvidente de la d\Žfinition de la fermeture initiale. Tout \Žl\Žment de $\udl I$ est $\leqs$ \ˆ  certain de $I$;
$$\text{ayant} \  I < \ovs I \ \text{on aura} \ \udl I < \ovs I \ \therefore \  \udl I \inq \ovs I\uds)$$
$$\text{de} \  F \inq \ovs I \ \text{vient}   \ \ovs I\uds) \inq \uds F$$

De ces relations r\Žsulte que les trous \Žtroits
$$\left[\ovs I \uds), \ovs I\right]\left[\uds F, \uds F\ovs)\right]$$
sont, dans le syst\me de trous de $(M <)$, incomparables puisque
$$\ovs I \uds) . \uds F\ovs) \inq \ovs I . \uds F \equiv 0$$
si la ligne $(L <)$ \Žtait compl\te. Quand ce ne serait pas le cas, et que l'un des ensembles $L_{if}$ dont on a d\Žj\ˆ parl\Ž pourrait \tre adjoint, le premier de ces trous pr\Žc\Žderait le second.

Chaque ligne compl\te de $(M <)$, commenc\Že en $I$, est prolong\Že par n'importe quelle autre  de $(\ovs I <)$.

En effet : soit $(L '<)$ une ligne qui passe par $I$, il est \Žvident que les \Žl\Žments de $L' - I$  v\Žrifient
$$I < L' - I \ \therefore \ L' - I \inq \ovs I$$

Soit  $(I <)$ un demi-rayon initial de $(M <)$, si $(F <)$ \Žtait une ligne compl\te de $(\ovs I  <)$, $(I + F <)$ le serait de $(M <)$.

En effet : que $(I + F <)$ est une ligne, c'est \Žvident. Qu'elle soit compl\te, parce qu'\ˆ l'int\Žrieur, ou \ˆ gauche de $(I <)$ aucun \Žl\Žment de $M - (I + F)$ ne rentre, ni \ˆ l'int\Žrieur ni \ˆ droite de $(F <)$. Ni entre $I$ et $F$ non plus, car il serait de $(\ovs I <)$, dans lequel la compl\te $(F <)$ a \Žt\Ž prise.

De ces propositions, r\Žsulte que
$$(I <) < (\ovs I <) \ \text{pour faire court} \ (I,\ovs I)$$
se pr\Žsente comme un croisement ascendant de $(I <)$ dans $(M <)$ puisque le chemin, commenc\Ž en $(I <)$, doit \tre poursuivi \ˆ l'int\Žrieur de $(\ovs I <)$ et peut \tre suivi par n'importe lequel des chemins
de cet ordre. Notez que $(I, \ovs I)$ ne sera g\Žn\Žralement pas un trou dans $(M <)$ parce que $(I < )$  sera  une section initiale seulement exceptionnellement.

Nous dirons {\it initialement identique} de deux lignes compl\tes, lorsqu'elles ont  un demi-rayon initial non vide commun. Dans le cas contrai\-re, nous les appellerons {\it initialement distinctes.}
Si deux lignes compl\tes sont initialement identiques, elles ont en commun un demi-rayon initial maximum $I$, qui est compl\Žt\Ž par deux demi-rayons finaux initialement distincts de $(\ovs I <)$.

En effet : soit $(L <)$ et $(R <)$ les deux lignes, et $I_j$ le demi-rayon initial g\Žn\Žrique commun. 
$$\sum I_j \equiv I$$
sera \Žgalement un demi-rayon initial commun; alors $I$ est l'un des $I_j$, et par construction, le plus grand. Par cons\Žquent
$$(L. \ovs I <) \ \et (R. \ovs I <)$$
seront initialement distincts dans $(\ovs I <)$.

{\it Indice cardinal} du croisement ascendant $(I, \ovs I)$, nous le dirons du nombre de lignes compl\tes de $(\ovs I <)$ initialement distincts deux \ˆ deux. 

De fa\c con analogue, en consid\Žrant un \Žl\Žment $a$ et en appelant $(a, \ovs a)$ son croisement ascendant,  on d\Žfinit l'indice cardinal du croisement ascendant de la m\me mani\re, \ˆ sa droite.

Si $(S <)$, totalement ordonn\Ž dans $(M <)$, \Žtait cofinal \ˆ $(I < )$,  demi-rayon initial, son prolongement ascendant serait form\Ž des lignes compl\tes de $(\ovs I <)$.

Effectivement : en raison de la suppos\Že cofinalit\Ž 
$$\udl S \equiv \udl I \ \therefore \ \ovs S \equiv \ovs I$$

\'Etant donn\Žs $(A, B)$ un trou de $(M <)$, et $(L <)$ une ligne compl\te, on v\Žrifie, en d\Žsignant par $N$ l'intervalle neutre de ce trou
$$L . A < L . N < L . B$$

En effet : d\Žsignant par $I_a$ l'\Žl\Žment g\Žn\Žrique de $L. A$, chaque \Žl\Žment $I$ de $L$ qui v\Žrifie $I <I_a$, sera de $A$; donc \Žgalement de $L. A$. L'\Žl\Žment $I_n$ de $L. N$, ne pouvant d\Žpasser aucun de $L. A$, v\Žrifiera
$$L.A < I_n \ \therefore \ L.A < L.N$$
Avec une d\Žmonstration, logiquement identique, on peut \Žtablir la seconde  affirmation.

Avec les hypoth\ses de la proposition ant\Žc\Ždente,
$$(L.A <) \ \text{est une ligne compl\te de} \ [(L . N)_* <]$$
$$\quad \qquad (L.N <) \ \text{est une ligne compl\te de} \ [(L . A)^*.(L. B)_* <]$$
$$(L.B <) \ \text{est une ligne compl\te de} \ [(L . N)^* <]$$

Effectivement :
$$\text{de} \   L. A < L.N \ \text{r\Žsulte} \ L. A \in (L . N)_*$$
et ainsi la premi\re proposition est v\Žrifi\Že, en y omettant, pour l'instant, le mot \guil complet\guir. Compl\te \ˆ l'int\Žrieur, et \ˆ gauche, car $(L <)$ est une ligne compl\te. Ses extensions \Žventuelles, avec des \Žl\Žments de $(L. N) _*$ occuperaient son extr\me droite, et donc, d\Žsignant par $x$ un tel \Žl\Žment,
$$L. A < x <L.N$$
et $(L <)$, contrairement \ˆ l'hypoth\se,  ne serait pas complet; alors un tel $x$ n'existe pas dans $(L. N)_*$. La troisi\me proposition, incluse dans l'\Žnonc\Ž, est corr\Žlative. Nous d\Žmontrerons la proposition interm\Ždiaire. De la proposition ant\Žc\Ždente, il suit
$$L.N \inq (L.A)^* . (L . B)_*$$
et \Žtant compl\te \ˆ l'int\Žrieur, si des \Žl\Žments de cet ensemble pouvaient y \tre ajout\Žs, ce serait faisable dans ses extr\Žmit\Žs seulement,  v. gr. [verbi gratia] $ x \in L. N$; mais 
$$\text{s'il est vrai que } \ x\in L. A \ \text{c'est que} \ L . A < x < L . N$$
et nous arriverions \ˆ la m\me contradiction que dans le cas pr\Žc\Ždent.

L'exemple que nous avons donn\Ž, pour montrer que la notion de ligne compl\te ne peut pas \tre
relativis\Že, prouve que $(L. A <)$  n'est pas forc\Žment compl\te dans $(A <)$.

Nous dirons qu'une ligne compl\te $(L <)$ {\it passe} par un trou $(A,B)$, lorsque  $L. N$ est vide; sinon qu'elle {\it ne passe pas}. 

\'Evidemment, lorsqu'une ligne compl\te passe par un trou $(A,B)$, elle se compose de, et tarit avec, les deux sections, initiale et finale, $L. A$ et $L. B$; il passera donc \Žgalement par $(\udl {L . A},\ovl {L .B)}$.

Il est clair que chaque ligne compl\te passe par chaque trou disjonctif.
C'est une condition n\Žcessaire et suffisante pour que la ligne compl\te $(L <)$ passe par le trou $(A, B)$, que l'ensemble $(L. A) ^*. (L. B)_*$  soit vide.

En effet :  $L. N$ \Žtant vide, $(L. A, L. B)$ sera un trou \Žtroit de    $(L <)$; mais, si $(L. A) ^*. (L. B) _*$ contenait un \Žl\Žment $m$ de l'ensemble $M$, on aurait 
$$L. A <m <L. B$$
et la ligne, contrairement \ˆ l'hypoth\se, ne serait pas compl\te; donc, la condition est n\Žcessaire. Supposons que l'ensemble $(L. N)$ ne soit pas vide : en vertu de ce qui a d\Žj\ˆ \Žt\Ž prouv\Ž, quelques lignes plus haut, le tout serait contenu dans $(L. A) ^*. (L. B) _*$  et ainsi la condition est suffisante.

\

13.--{\sc Transversales compl\tes.}

\

Les  sous-ensembles de $M$, dont les paires sont toutes incomparables, nous les appellerons {\it transversales}. La structure $(M <)$, dans laquelle elles se trouvent, n'agit pas sur elles,  les  laisse dans un \Žtat amorphe, en relation avec leur environnement. On dira que la  transversale $T$  est compl\te lorsque tout \Žl\Žment de $M - T$ est comparable \ˆ l'un de $T$. L'ensemble $\uds t + \ovs t$ comprenant tous les \Žl\Žments comparables \ˆ $t$, l'\Žquation
qui d\Žfinit la compl\tude de $T$ sera 
$$T \equiv M - \sum_{t\in T} (\uds t + \ovs t)$$

L'ensemble de tous les \Žl\Žments maximaux de $(M <)$ constitue
une transversale.

En effet :  $\ovs a$ et $\ovs b$ \Žtant vides, lorsque $a$ et $b$ sont deux maxima de $(M <)$, aucune des deux relations suivantes n'est vraie
$$a <b \  \text{ni} \ b < a$$

Cette transversale n'est pas n\Žcessairement compl\te, et ainsi, \ˆ proprement parler, il se peut que l'on ait
$$\udl T \inc M$$
et, dans ce cas, $(M - \udl T <)$ n'a pas de maximum (*). C'est ce qui se passe dans l'ordre d\Žfini dans la figure 4.

Soit $T$ une transversale compl\te  et $z$  l'un des \Žl\Žments de $M - T$, l'un des ensembles $(T. \uds z),  (T. \ovs z)$ a des \Žl\Žments et l'autre pas.

En effet : que les deux aient des \Žl\Žments est impossible car on aurait
$$x < z < y \ \therefore \ x < y$$
  tous deux \Žtant  de $T$. Si les deux \Žtaient vides, $z$ serait incomparable \ˆ tous ceux de $T$, et celle-ci ne serait pas compl\te, contrairement \ˆ  hypoth\se. (**)

\
-----------------

(*) Prenant la transversale des minima, que $M - \ovl T$ soit vide ou non, pour la structure, partiellement ordonn\Že par l'inclusion, des bases d'une science math\Žmatique, provenait la classification de la p.108, de notre article \guil Estructuras deductivas \guir[Les structures d\Žductives]
 (Rev. Mat. Hisp-Am. (IV) {\it 14} (1954) (104-117).

(**) Le lecteur notera l'analogie avec les principes de contradiction et du t. exclu [le tiers exclu].

-------------------

 \
 
 En vertu de la proposition pr\Žc\Ždente, chaque  transversale compl\te $T$, permet de r\Žpartir les
\Žl\Žments de $M$ en trois classes : 

\noi {\it Classe} 0; ceux de $T$.
 
\noi {\it Classe} 1; ceux de $\udl T - T$, qui constituent une section initiale.

\noi {\it Classe} 2; ceux de $\ovl T - T$, qui constituent une section finale. 

Il convient de noter qu'en g\Žn\Žral,
$$(\udl T -T, \ovl T - T)$$
ne constitue pas un trou de $(M <)$ car il n'est en g\Žn\Žral pas vrai que
$$\udl T -T < \ovl T - T$$
ainsi, dans la figure 8, pour aucune transversale ces ensembles ne constituent de trou, comme on peut le voir
facilement.

Dans $(\udl T <)$ les \Žl\Žments de $T$ constituent la transversale compl\te  de leurs
maxima. Il ne faut pas en conclure que chaque ligne compl\te de $(\udl T < )$ se termine par un \Žl\Žment de $T$ car ellle pourrait avoir la structure ordinale de la figure 9, dont la transversale des maxima, ici compl\te, est constitu\Že par les points $a_i$, $i$ parcourant les nombres entiers. L'horizontale de cette figure constitue une ligne compl\te, qui ne se termine par aucune desdites transversales.

Une transversale compl\te r\Žpartit les \Žl\Žments d'une ligne compl\te $(L <)$ en trois classes
$$(\udl T - T) . L \ , \  T. L \ , \ (\ovl T-T). L$$
la premi\re, avec l'ordre subordonn\Ž, constitue une section initiale,

\

\includegraphics[width=4.8in]{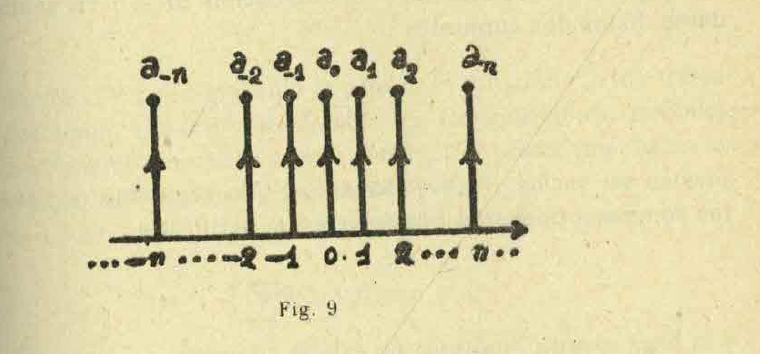}

\

\noi la derni\re de mani\re analogue  une finale; l'interm\Ždiaire a au plus un \Žl\Žment et pourrait \tre vide. Sur la ligne $(L <)$
$$[(\udl T -T). L, (\ovl T - T) . L]$$
constitue un trou, qui est disjonctif lorsque 
$$T . L \equiv  {\rm O}$$
En ce sens, on peut dire qu'une transversale traverse une ligne compl\te passant par un point, celui de $T. L$ quand il y en a un, ou par le trou disjonctif, d\Žcrit ci-dessus, si $T. L$ \Žtait vide. 

Si la transversale compl\te $T$  croisait la ligne compl\te $(L <)$ par son \Žtroit $(L', L'')$,  $l'$ et $l''$   \Žtant respectivement les \Žl\Žments g\Žn\Žriques de $L'$ et $L''$,
$$l' < x   \et  y < l''$$
ont des solutions en $T$, qui constituent les ensembles
$$T. \ovs {l'}  \et  T . \uds l''$$
En parcourant $I'$, dans le sens croissant, $(L' <)$,  il y aura certaines fois, d'autres non, l'ensemble limite d'\Žl\Žments communs \ˆ tous, ceux de l'ensemble 
$$T . \ovs {L'} \equiv \prod T . \ovs{l'}$$
De mani\re analogue, en parcourant $I''$  dans le sens d\Žcroissant  de l'ordre $(L'' <)$. Ces deux ensembles
$$T. \ovs{L'} \ \et T. \uds{L''}$$
peuvent \tre vides ou poss\Žder des \Žl\Žments. Ils n'ont s\žrement pas d'\Žl\Žments communs car celui qui renfermerait $x$ v\Žrifierait
$$L' < x < L''$$
et la ligne, contrairement \ˆ l'hypoth\se, ne serait pas compl\te.

Prenant dans $(\ovs T <)$ une transversale compl\te  et l'\Žtendant jusqu'\ˆ ce que, \ˆ l'int\Žrieur de $(M <)$, elle soit aussi compl\te, on obtiendrait une transversale compl\te $T'$. Serait vraie 
$$T < T'$$
puisque,  $t$ et $t'$ \Žtant de leurs \Žl\Žments respectifs, s'il se trouvait  parfois
$$t' \leqs t$$
par la d\Žfinition de $T'$, il y aurait des \Žl\Žments post\Žrieurs \ˆ tous ceux de $T$; on verrait ainsi qu'il y avait des \Žl\Žments comparables dans $T'$.

{\it Probl\me}. -- Existe-t-il des ordres dans lesquels $L. T$, quels que soient $L$ et $T$, avec la signification pr\Žc\Ždente, poss\de  un \Žl\Žment ?

Bien s\žr, que $L . T$ soit toujours vide est impossible car, $t$ \Žtant un \Žl\Žment de $T$, il y a des lignes qui le traversent.

\newpage

14.--{\sc Bases lin\Žaires d'un ordre.}

\

Les figures suggestives utilis\Žes dans cet article indiquent l'importance de la formation g\Žn\Žrale des ordres, non seulement en commen\c cant par d'autres totaux \ˆ deux \Žl\Žments, mais avec n'importe quel nombre, en particulier avec ceux qui sont des lignes compl\tes de la structure.

Si $L_i$, lorsque $i$ parcourt un ensemble convenable $I$, nous fournissait la totalit\Ž des lignes compl\tes de $(M <)$, il est clair que, par la confusion des ordres du syst\me
$$(L_i <) \ \pour i\in I$$
nous obtiendrions $(M <)$. Ce qui est int\Žressant est que, puisque la loi transitive constitue un m\Žcanisme permettant de {\it d\Žduire} de nouvelles relations d'ordre, de la g\Žn\Žration pr\Žc\Ždente, toutes les lignes compl\tes ne sont g\Žn\Žralement pas utilis\Žes, mais seule une partie d'entre elles suffit.

Appelons $J$ le sous-ensemble g\Žn\Žrique de $I$  : les ordres 
$$\left(\sum_{j \in J} L_j < \right) \ J \inq I$$
seront en g\Žn\Žral plus faibles que $(M <)$; lorsqu'il co\•ncide avec elle, nous dirons que le syst\me de lignes compl\tes
$$(L_j  <) \ j\in  J$$
est une {\it base lin\Žaire} de l'ordre $(M <)$. Nous d\Žsignerons dans la suite de mani\re g\Žn\Žrique l'ensemble $J$ par $B$.

En ordonnant le syst\me $I_1$ des sous-ensembles de $I$, par inclusion stricte, le syst\me $\beta$ des bases lin\Žaires, constitue une section finale de $(I_1 \inc)$, et, lors du parcourt  de la structure $(I_1 \inc)$, dans le sens croissant, qui finit par entrer dans $\beta$ car le m\me ensemble $I$, qui est le plus grand de $(I_1 \inc)$, constitue une base lin\Žaire. En particulier, en suivant les lignes
compl\tes de $(I_1 \inc)$, elle finit par entrer dans la section finale des bases lin\Žaires. 

Il peut arriver que $(\beta \inc)$ ait un \Žl\Žment infime; l'ensemble qui le  constitue sera une {\it base absolue}  de l'ordre $(M <)$. Si cela ne se produit pas,  $(\beta \inc)$  peut encore avoir  des \Žl\Žments minimaux; nous les appellerons {\it bases irreductibles} de l'ordre $(M < )$.
Bien s\žr, a priori, celles-ci pourraient \Žgalement manquer, et alors chaque base en contiendrait une autre de l'ordre $(M <)$. Ayant des bases irr\Žductibles, deux cas peuvent \tre distingu\Žs : {\it a}) les minima de $(\beta \inc)$ constituent une transversale compl\te de $(\beta \inc)$, et, par cons\Žquent, toute base lin\Žaire de $(M <)$ contient une base irr\Žductible; {\it b}) la transversale des bases irr\Žductibles est incompl\te dans $(\beta \inc)$, et par suite, en plus des bases lin\Žaires qui contiennent des bases irr\Žductibles, il en existe d'autres qui ne contiennent aucune des autres bases irr\Žduc\-tibles, sinon elles-m\mes, et toutes les autres qui y sont contenues, sont r\Žductibles.

\

\includegraphics[width=4.8in]{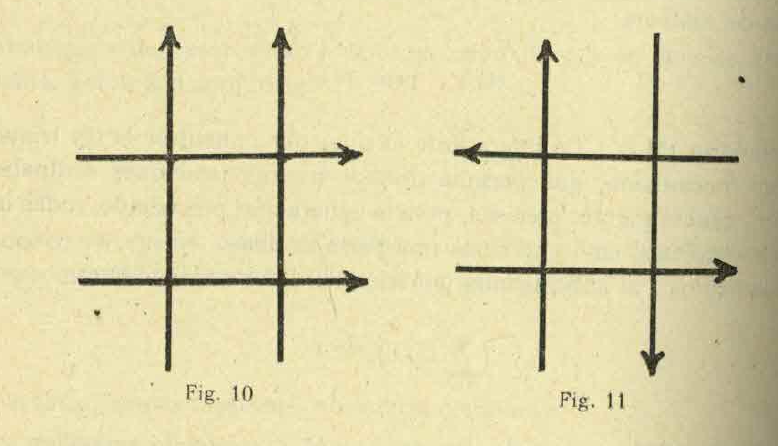}

\

Si tous ces cas sont effectivement r\Žalisables, est une question probl\Ž\-matique.

Les exemples suivants clarifient les notions ant\Žrieures. Les figures 3, 10, 11, nous montrent des ordres, tous engendr\Žs  par quatre lignes compl\tes, dont le type ordinal est celui de l'ordre r\Žel. 
Le nombre total de lignes compl\tes est, cependant, diff\Žrent; dix-neuf pour la figure 3, dix-huit pour la figure 10, seize pour la figure 11. C'est suffisant pour  voir l'anisomorphisme de ces ordres. Les transversales, dont le nombre est infini dans chacune d'elles, se compose au plus de quatre \Žl\Žments, mais il y en a avec trois et avec deux; aucun dans la figure 3, un dans la figure 10, deux dans la figure 11, ce qui prouverait \Žgalement l'anisomorphisme desdits ordres. Les ordres d\Žfinis par les figures 12 et 13, avec le premier $a$, et le dernier $b$, sont isomorphes pour l'ordre, m\me si topologiquement, et sans sortir du plan, on ne peut les superposer.

Au lieu d'une fusion conjointe, qui suppose des \Žl\Žments communs dans les ordres qui se confondent, on pourrait partir d'un syst\me d'ordres $(B_i <)$ -- totaux dans le cas int\Žressant que nous examinons maintenant -- dont les ensembles n'ont pas  d'\Žl\Žments communs, deux \ˆ deux. Ceci serait r\Žalis\Ž en d\Žfinissant, de toutes les mani\res possibles, une structure de classification, pour l'ensemble r\Žunion $\sum B_i$, qui n'agit pas librement, car on doit \Žviter les cycles, ce qui n\Žcessite de la limiter par les conditions suivantes :

\

\includegraphics[width=4.8in]{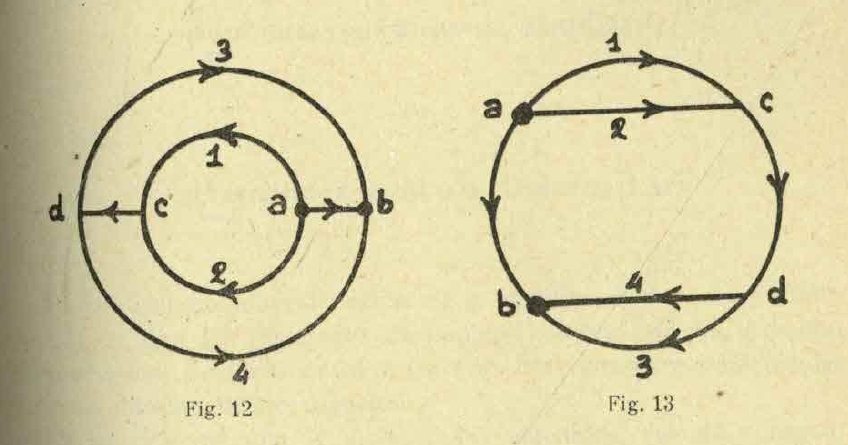}

\

1.$^a$) L'intersection de chaque classe avec chacun des ensembles $B_i$ d\Žj\ˆ ordonn\Žs  n'aura jamais deux \Žl\Žments ou plus.
 
2.$^a$) Les classes $P$ et $Q$ , \ˆ leur intersection avec $B_i$ fourniront deux \Žl\Žments $p_i$ et $q_i,$ dont la relation d'ordre reste invariante en faisant varier $i$ de toutes les mani\res possibles.

Ce proc\Žd\Ž nous donnerait tous les ordres possibles, dont la base est constitu\Že des ordres $(B_i <)$ du syst\me donn\Ž.

\

\

\centerline{\bf  FIN DE LA TRADUCTION}

\

\

\centerline{\bf Table des mati\res}

\head{Prologue\hfill{6}}

\head{CHAPITRE 1 LES STRUCTURES BINAIRES \hfill{7}}
				
\

\noi1.--{\sc Les structures binaires en g\Žn\Žral}\hfill{7}

\noi 2.--{\sc Les structures binaires transitives} \hfill{10}

\noi 3.--{\sc La classification des structures binaires transitives} \hfill{14}

\noi 4.--{\sc Ordre et classification engendr\Žs par une structure $(M \leqs)$.} \hfill{15}

\head{CHAPITRE II. LES ORDRES \hfill{17}} 

\noi 5.--{\sc Notions pr\Žliminaires.} \hfill{17}

\noi 6.--{\sc Majorant et minorant d'un ensemble.}\hfill{20}

\noi 7 .--{\sc Fermetures initiale et finale d'un ensemble.} \hfill{21}

\noi 8.--{\sc Extension syst\Žmatique des ordres.}\hfill{24}

\noi 9.--{\sc Les diff\Žrents types de trous.} \hfill{28} 

\noi 10.--{\sc Ordre naturel des trous.} \hfill{29}

\noi 11.--{\sc Matrice d'ordre pour un nombre cardinal.} \hfill{32}

\noi 12.--{\sc Lignes compl\tes.} \hfill{34}

\noi 13.--{\sc Transversales compl\tes.} \hfill{38}

\noi 14.--{\sc Bases lin\Žaires d'un ordre.} \hfill{42}

\

\ 

\centerline{\bf VOCABULAIRE}

\

\noi Dans m p.x, le m renvoie \ˆ la page de la traduction, le 
p.x \ˆ la page du texte original

\

r\Žflexif 7 p.104

irr\Žflexif 7 p.104

sym\Žtrique 7 p.104

asym\Žtrique 7 p.104

incomparable 7 p.104

comparable 7 p.105

plus faible 7 p.105

non connect\Ž 8 p.105

connexe 8 p.105

r\Žunir 9 p.105

sousconjoindre 9 p.106

ligne compl\te 9 p.107

transversale compl\te 9 p.107

fermer \ˆ gauche 9 p.107

inverser 10 p.107

noyau structur\Ž 10 p.107

r\Žsidu amorphe 10 p.107

isomorphe au sens large 10 p.107

proposition 11 p.108

vraies 11 p.108

fausse 11 p.108

crit\re de v\Žrit\Ž 11 p.108

m\Žcanisme d\Žductif 11 p.108

fermeture d\Žductive 12 p.109

base 12  p.110

r\Žductible 12 p.110

irr\Žductible 12 p.110

base absolue 12 p.110

con-fusion 13 p.111

odre total 18 p.115

ordre partiel 18 p.115

maximum 18 p.116

minimum 18 p.116

supremum 18 p.116

infimum 18 p.116

inconfusible 19 p.117

cycle 19 p.117

finalement sup\Žrieur 21 p.119

fermer initialement 21  p.120

fermeture finale 21 p.120

section initiale 22 p.120

section finale 22 p.120

intervalle 22 p.120

ordre ramifi\Ž 23 p.121

cofinaux 23 p.121

co\•nitiaux 23 p.121

enveloppe sup\Žrieurement 23 p.122

extension 24 p.123

sous-ordre 25 p.123

sur-ordre 25 p.123

extension limite 25 p.124

intervalle neutre 27 p.127

trou ext\Žrieur 27 p.127

trou couvert 27 p.127

trou appuy\Ž 27 [estribado] p.127

trou \Žtroit 27 p.127

trou disjonctif 27 p.128

ligne 34 p.134

demi-rayon initial 35 p.135

demi-rayon final 35 p.135

indice cardinal 37 p.137

passe par un trou 38 p.138

ne passe pas 38 p.138

transversale 38 p.139

base lin\Žaire 42 p.143

base absolue 42 p.143

base irr\Žductible 42 p.144

\

\

\

\

\head{Quelques pr\Žcisions}

\

Apr\s les deux points : on commence, en g\Žn\Žral, par une minuscule comme la  typographie du fran\c cais l'exige.

On a essay\Ž de s'en tenir, autant que possible,  aux traductions suivantes.

pues : car

por tanto : donc

tendr\'ase : on aura

al pasar : en passant

desde luego : bien s\žr

por tener : en ayant

luego : alors

o sea : c'est-\ˆ-dire

siendo : \Žtant ou soit

as\'i : ainsi

recorrere : parcourir

acierto : convenable (?)

agotar : tarir

constar con : se compose de

\

\

Le style de {\sc Cuesta} nous a sembl\Ž particulier. Aux fins de comparaison avec des math\Žmaticiens contemporains voici un
extrait, la page 77, du livre de {\sc Fernando Hern\'andez Hern\'andez}, [{\sl Teor\'ia de Conjuntos (una introducci\'on)},  Universidad Nacional Autonoma de Mexico, 1998],   suivies de deux pages, 104 et 105, de {\sc Cuesta}.

\includepdf{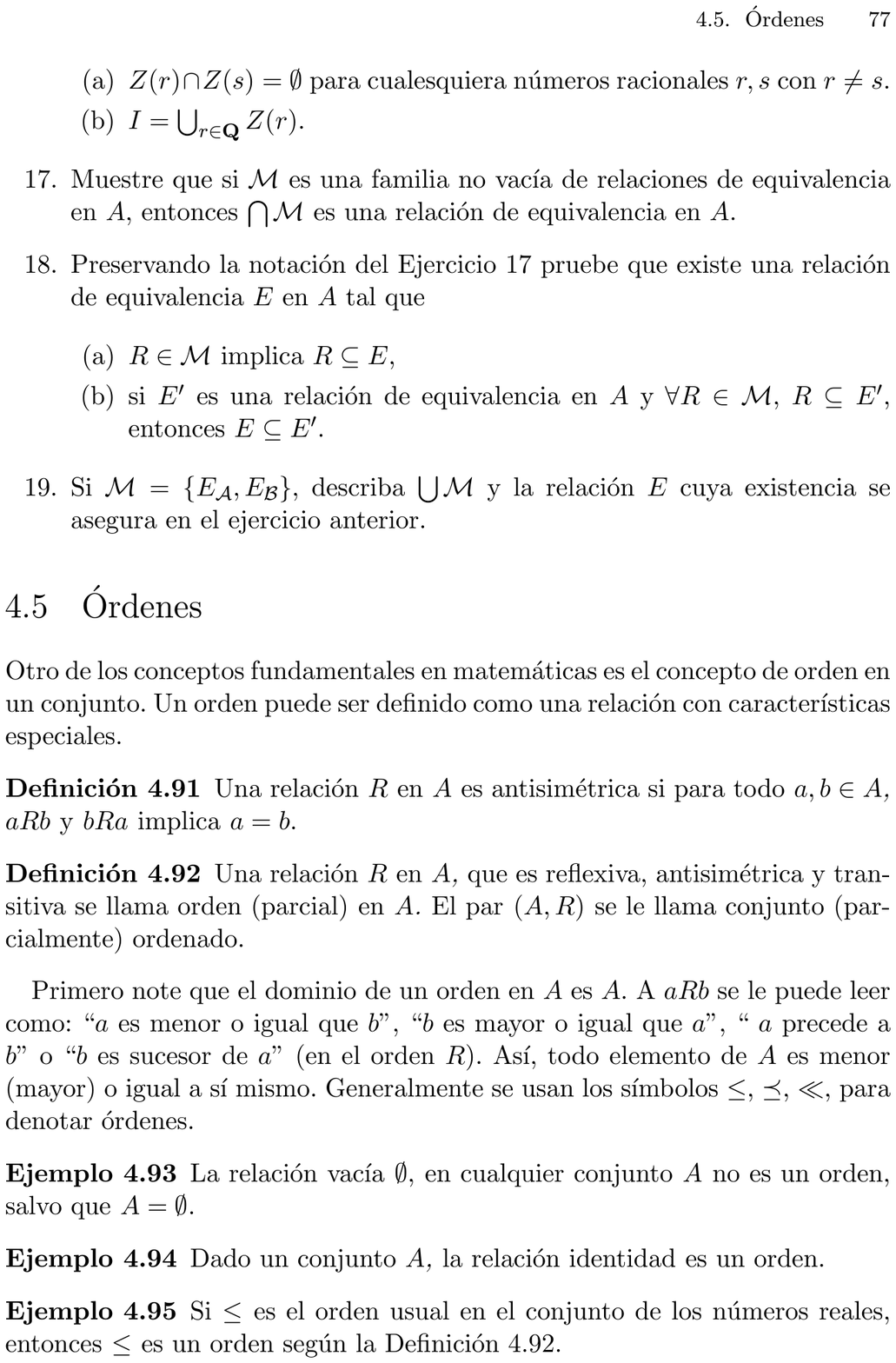}

\includepdf{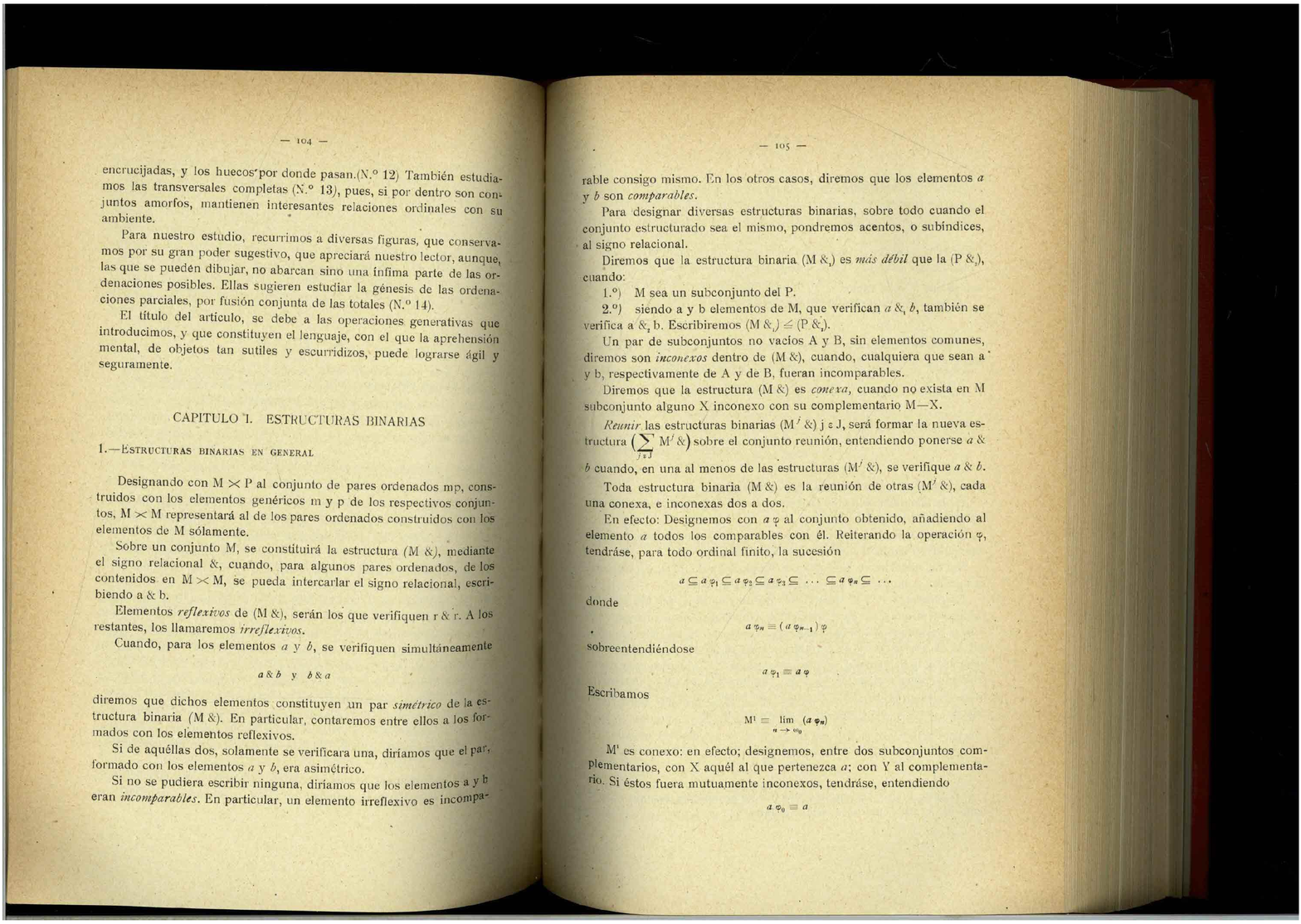}

\

\

\
 
\

\
\enddocument